\documentclass[ijoc,nonblindrev]{StyleFiles/informs3} 

\OneAndAHalfSpacedXI

\usepackage[font=small]{caption}

\usepackage{StyleFiles/custom_formats}
\usepackage{StyleFiles/custom_ian}
\theoremstyle{definition}

\usepackage[skip=15pt]{subcaption}

\usepackage[noend]{algpseudocode}
\usepackage{etoolbox}

\newcounter{subroutine}


\usepackage{natbib}
 \bibpunct[, ]{(}{)}{,}{a}{}{,}%
 \def\bibfont{\small}%
\TheoremsNumberedThrough     

\EquationsNumberedThrough    

\begin{document}

 \RUNAUTHOR{Bodur, Chan, and Zhu}

\RUNTITLE{Network Flow Models for Adaptive Robust Optimization}

\TITLE{Network Flow Models for Robust Binary Optimization with Selective Adaptability}

\ARTICLEAUTHORS{%
\AUTHOR{Merve Bodur}
\AFF{School of Mathematics, University of Edinburgh, Edinburgh  EH9 3FD, UK,
\EMAIL{merve.bodur@ed.ac.uk}}
\AUTHOR{Timothy C. Y. Chan}
\AFF{Department of Mechanical and Industrial Engineering, University of Toronto, Toronto, Ontario M5S 3G8, Canada, 
  \EMAIL{tcychan@mie.utoronto.ca}
  }
\AUTHOR{Ian Yihang Zhu}
\AFF{NUS Business School, National University of Singapore, Singapore 119245,
\EMAIL{ianyzhu@nus.edu.sg}}
}

\ABSTRACT{Adaptive robust optimization problems have received significant attention in recent years, but remain notoriously difficult to solve when recourse decisions are discrete in nature. In this paper, we propose new reformulation techniques for adaptive robust binary optimization (ARBO) problems with objective uncertainty. Without loss of generality, we focus on ARBO problems with ``selective adaptability", a term we coin to describe a common class of linking constraints between first-stage and second-stage solutions. Our main contribution revolves around a collection of exact and approximate network flow reformulations for the ARBO problem, which we develop by building upon ideas from the decision diagram literature. Our proposed models can generate feasible solutions, primal bounds and dual bounds, while their size and approximation quality can be precisely controlled through user-specified parameters. Furthermore, and in contrast with existing solution methods, these models are easy to implement and can be solved directly with standard off-the-shelf solvers. Through an extensive set of computational experiments, we show that our models can generate high-quality solutions and dual bounds in significantly less time than popular benchmark methods, often by orders of magnitude.}


\maketitle

\vspace{-0.7cm}

\section{Introduction}

Robust optimization (RO) has become a well-established paradigm for modeling and solving decision-making problems under uncertainty. It has found diverse applications across a wide range of problem domains, and is particularly well-suited for environments where the distribution of parameters in a model are difficult to characterize or when there is a need to consider worst-case outcomes. While early research has predominantly focused on static RO models, where a single robust solution is generated for all possible parameter realizations \citep{bertsimas2013robust}, there has been significant recent interest in \emph{adaptive} robust optimization, where \emph{recourse} decisions can be made once additional information is revealed \citep{yanikouglu2019survey}.

Adaptive robust optimization problems, particularly those with discrete recourse variables, have numerous applications. Common examples include, but are not limited to, routing \citep{eufinger2020robust}, scheduling \citep{yan2018robust}, facility location \citep{hanasusanto2015k}, network design \citep{alvarez2015recoverable}, batching \citep{bayram2022optimal}, assignment \citep{dacs2020review} and matching problems \citep{mcelfresh2019scalable}. Despite their importance, adaptive robust optimization models can be significantly more challenging to solve relative to their static counterpart. While static models can typically be reformulated and solved as tractable mixed-integer linear programming (MILP) problems, these reformulation techniques do not extend to adaptive models, where optimal discrete recourse decisions must be defined for each and every parameter value lying within an uncertainty set. 
Few solution methods exist for solving these models, and most require elaborate and carefully-tuned iterative algorithms \citep{zeng2013solving, kammerling2020oracle, arslan2022decomposition}. Significant attention has instead been placed on approximation techniques, particularly ones that yield MILP formulations as these can be easily implemented through standard commercial solvers \citep[e.g.,][]{hanasusanto2015k}.



In this paper, we consider adaptive robust binary optimization (ARBO) problems with objective uncertainty. We focus, without loss of generality, on ARBO problems with \emph{selective adaptability}, a term that we coin to describe a class of linking constraints where first-stage variables fix the values of some recourse variables without restricting the rest; that is, the remaining recourse variables have the ability to adapt to new information. This distinct structure underlies a variety of planning and sequential decision-making tasks, 
and we highlight several examples in Section \ref{sec:model_intro}. More importantly, the particular structure of these linking constraints motivates the design of new reformulation techniques and solution approaches. In particular, we employ ideas from the decision diagram community \citep{castro2022decision} to reformulate these ARBO problems into various exact and approximate \emph{constrained network flow models}. These models are flexible, easy to implement, can be solved directly using standard commercial solvers, and offer a number of computational advantages compared to existing methods. 

While many common problems naturally fit the description of ARBO with selective adaptability, we emphasize that our focus on these problems comes without loss of generality. Specifically, we will show that any ARBO problem can be reformulated into one that has this property through the use of auxiliary variables. Thus, all reformulation techniques and models presented in this paper are relevant to the general class of adaptive robust binary problems with objective uncertainty.


A concise list of our main contributions are as follows:

\begin{enumerate}
    \item We introduce ARBO models with selective adaptability, and show how the structure of the linking constraints can be exploited when we have a description of the convex hull of the recourse feasible space. By highlighting the connection between decision diagrams and this convex hull description, we show that the ARBO models can be reformulated as constrained network flow problems, where recourse decisions are represented by continuous flow along certain links in a large capacitated network, and where capacity constraints are given by the values of first-stage decisions. We show that these models have MILP formulations and can thus be directly solved using standard commercial solvers. 
    \item We introduce three approximation techniques to generate smaller and more compact network flow models for large ARBO problems. The first two techniques utilize approximate decision diagrams to formulate inner and outer approximations of the recourse feasible space, and thus generate primal and dual bounds for a given ARBO problem, respectively. The third technique, which we term as a generalized multi-network flow model, pools together a collection of multiple constraints and networks to generate a dual bound. We also outline new methods for specifying the size and quality of these approximation models. 
    \item We examine the performance of our exact and approximate network flow formulations in two sets of numerical experiments spanning a robust project investment problem and a robust assignment problem. For smaller ARBO problems, we find that the exact flow formulations are small in size and can be solved efficiently. For larger ARBO problems, we show that approximate network flow formulations are substantially smaller and easier to generate. More importantly, these approximate formulations consistently generate near-optimal solutions and high-quality dual bounds in solution times that are orders of magnitude lower than exact formulations and benchmark models. Furthermore, by adjusting the size of the approximate models, we show that the complexity and solution times of the approximate models can be reduced substantially while sacrificing little in terms of solution quality, making these methods highly scalable for large ARBO problems.
\end{enumerate}

The rest of the paper unfolds as follows. In Section \ref{sec:lit_review}, we review the relevant literature on robust optimization, decision diagrams, and network flow models. In Section \ref{sec:model_intro}, we introduce the ARBO model and the concept of selective adaptability. In Section \ref{sec:BDD}, we present the exact network flow reformulation, while Section \ref{sec:approximations} presents a series of approximation techniques that result in more compact and tractable network flow formulations. Section \ref{sec:experiments} examines the performance of the network flow models across various numerical experiments. We conclude in Section \ref{sec:conclusion}.

We summarize key notation used in this paper. All vectors and matrices are in bold letters, while sets are denoted using calligraphic letters. Let $\mS$ be a feasible set described using linear and integrality constraints. We use $\rel(\mS)$ to describe the relaxation of $\mS$ obtained by removing integrality constraints, and $\conv(\mS)$ to denote the convex hull of $\mS$; by definition, $\conv(\mS) \subseteq \rel(\mS)$. The set of extreme points of any polyhedral set $\mS$ is denoted by $\ext(\mS)$.

\section{Literature Review}\label{sec:lit_review}

We first review the relevant literature for robust binary optimization with objective uncertainty. We then review the relevant literature on decision diagrams and network flow models. We refer readers interested in a general overview of robust optimization to the comprehensive surveys by \cite{bertsimas2013robust} and \cite{yanikouglu2019survey}. 

\subsection{Adaptive robust binary optimization}\label{subsec:RO_litreview}

Static robust binary optimization (RBO) problems with objective uncertainty are prevalent and have been well-studied in the literature, and we refer to \cite{kasperski2017robust} and \cite{buchheim2018robust} for comprehensive surveys of these problems. To summarize, RBO problems are challenging, and even robust formulations of polynomially-solvable combinatorial optimization problems can be NP-hard \citep{buchheim2018robust}. Nonetheless, they can be formulated and solved as MILP models by reformulating the inner adversarial problem using standard optimality conditions. However, these reformulations do not extend to the adaptive setting, where recourse decisions must be defined for each and every parameter vector within the uncertainty set.

Few algorithms exist for solving general two-stage RBO problems with objective uncertainty. The first and most well-known approach is the nested constraint-and-column generation method proposed by \cite{zeng2013solving}. This algorithm iterates between a master problem and a subproblem, the latter iteratively adding variables and constraints into the master problem. The subproblem is itself a max-min problem, and is commonly solved using a decomposition algorithm (e.g., constraint-and-column generation) for each outer iteration. A second approach is the branch-and-price algorithm developed by \cite{arslan2022decomposition}. Within each node of a branch-and-bound search tree, the algorithm solves a series of iterative pricing problems to add columns that correspond to feasible second-stage solutions for a fixed first-stage solution generated by a master problem defined within each given node. Finally, \cite{kammerling2020oracle} define a specialized branch-and-bound algorithm that branches over first-stage solutions while using another algorithm to iteratively refine the dual bound in each node of the branch-and-bound tree. The first approach of nested constraint-and-column generation is known to scale poorly to larger instances, because each iteration of the algorithm adds a full copy of all variables and constraints of the recourse problem to the master problem \citep[e.g.,][]{dumouchelle2023neur2ro}. The algorithm also may not converge in a finite number of iterations, because binding constraints may exist in the interior of the uncertainty set \citep{arslan2022decomposition}. On the other hand, the latter two approaches are iterative algorithms that require carefully-tuned intermediary steps and may be confined to specific software; for example, many commercial solvers limit the degree of user customization of specific branch-and-bound processes.

The challenges that surround these exact solution methods have motivated much interest in approximation techniques (e.g., \cite{vayanos2011decision, bertsimas2015design, postek2016multistage}; see \cite{arslan2022decomposition} for a comprehensive review). For example, \cite{bertsimas2015design} and \cite{bertsimas2016multistage} consider the use of piecewise constant decision rules to represent the binary recourse variables. On the other hand, \cite{hanasusanto2015k} and \cite{subramanyam2020k} propose the $K$-adaptability approximation method which pre-identifies $K$ recourse solutions to implement in the second stage. The $K$-adaptability model, which can be solved as a MILP model, implicitly partitions the uncertainty set into a finite number of subsets for which a single recourse decision is assigned to each. In recent years, this MILP formulation has become a commonly used approximation method, particularly for the problem class that we will examine (e.g., \cite{dumouchelle2023neur2ro, arslan2022decomposition}). 


In this paper, we present a new set of exact and approximate MILP reformulations for ARBO problems with objective uncertainty. These formulations, which are presented in the form of constrained network flow problems, are flexible, easy to formulate and implement in standard solvers, and generate high-quality solutions and bounds in times that can be significant lower than alternative solution methods.

\subsection{Decision diagrams and network flow models}

Decision diagrams (DDs) are increasingly used to generate solutions and bounds for discrete optimization problems. With few exceptions, this literature focuses on deterministic problems \citep{bergman2016decision}. We refer to \cite{castro2022decision} for a comprehensive survey of this literature.

Several works have examined the use of decision diagrams for two-stage stochastic programming with binary recourse decisions. \cite{lozano2018binary}, \cite{guo2021logic}, and \cite{macneil2022leveraging} consider various two-stage problems for which the recourse feasible space under each scenario can be represented as a decision diagram, making the two-stage problem amenable to standard Benders decomposition techniques that would otherwise only apply to problems with continuous recourse  \citep{rahmaniani2017benders}. In their problems, first-stage decision variables serve as capacity constraints in the diagrams, and given a first-stage decision, a shortest path algorithm can be applied on the diagram to generate a cut for the master problem. On the other hand, \cite{serra2019last} consider a monolithic formulation of a two-stage scheduling problem where a subset of constraints are replaced with DD-based network flow constraints. Integrality of the flows across these networks are then enforced using binary variables. 

To our knowledge, there is only one application of decision diagrams for robust optimization. \cite{lozano2022constrained} consider the reformulation of static robust binary optimization models as DD-based constrained shortest path problems for which specialized algorithms can be employed. 

In contrast to these works, we consider the first use of DD-based network flow models for adaptive robust optimization. We show that in our problem setting, the complete set of optimal recourse decisions given any first-stage solution can be represented as a continuous flow across a network (i.e., rather than a flow along a single path). We extend the analysis by proposing several primal and dual bounding techniques using the concept of approximate decision diagrams \citep{castro2022decision}. These are the first application of DD-based approximation schemes for robust optimization, and we show that they can lead to highly tractable and effective models in our numerical experiments.

Finally, we remark that decision diagrams are also closely related to the state transition graphs found in the dynamic programming literature; we refer to \cite{hooker2013decision} and \cite{castro2022decision} for detailed comparisons of these two research areas. \cite{de2022arc} provide a survey on dynamic programming-based network flow formulations for deterministic optimization problems. As the authors point out, one of the main advantages of these formulations is that they can be solved directly using standard MILP solvers, overcoming the necessity of more elaborate iterative methods. We consider exact and approximate DD-based network flow formulations for similar reasons, and show that they are also highly effective for our robust optimization problems.

\section{Robust Binary Optimization with Selective Adaptability}\label{sec:model_intro}

In this section, we define the ARBO problem of interest, introduce the concept of selective adaptability, and derive several preliminary insights based on the structure of the ARBO problem.

\subsection{Problem definition}\label{subsec:prelim}

Consider an ARBO problem of the form
\begin{align}\label{model:1}
\underset{\bx \in \mX}{\text{min}} \ \, \underset{\bxi \in \Xi}{\text{max}} \ \, \underset{\by \in \mY \cap \mS(\bx)}{\text{min}} \ \, \bc^\top \bx + \bxi^\top \by,
\end{align}
where $\mX \subseteq \{0,1\}^m$ defines the feasible set of binary first-stage decisions, $\Xi := \{\bxi \, | \, \bT \bxi \leq \bd\}$ is a bounded polyhedral uncertainty set, 
and $\mY \cap \mS(\bx) \subseteq \{0,1\}^n$ defines the feasible set of  binary recourse decisions where  
the set $\mY := \{\by \in \{0,1\}^n \, | \, \bG \by \geq \bh \}$ captures the  constraints that do not depend on $\bx$, while $\mS(\bx) \subseteq \mathbb{R}^n$ models the linking constraints. In this paper, we focus on problems where $\mS(\bx)$ satisfies a condition that we term \emph{selective adaptability}.

\begin{defn}[Selective Adaptability]
A linking constraint between first-stage and second-stage decision variables is defined to be selectively adaptive if it can be expressed as $y_i \leq x_j$, $y_i = x_j$ or $y_i \geq x_j$, for some $i \in \{1, \ldots, n\}$ and $j \in \{1, \ldots, m\}$. The ARBO problem  \eqref{model:1} is said to have selective adaptability if all constraints defining $\mS(\bx)$ are selectively adaptive, i.e., 
\begin{align}\label{cons:linking}
    \mS(\bx) = \left\{\by \in \mathbb{R}^n \ \middle | \begin{array}{l}
    y_{i} \leq x_j \hspace{1.2cm} \forall (i,j) \in \mU_1\\
    y_{i} = x_j \hspace{1.2cm} \forall (i,j) \in \mU_2\\
    y_{i} \geq x_j \hspace{1.2cm} \forall (i,j) \in \mU_3\\
    \end{array} \right\},
\end{align}
where $\mU_1, \, \mU_2, \, \mU_3$ are disjoint sets, and $\mU_1 \cup \mU_2 \cup \mU_3 \subseteq \{1,\ldots,n\} \times \{1, \ldots, m\}$.
\end{defn}

An ARBO problem with selective adaptability implies that the value of each binary recourse variable $y_i \in \{0,1\}$ is either (i) fixed by the values of first-stage decisions $\bx \in \mX$ (due to constraints corresponding to $\mU_2$, and those associated with $\mU_1$ and $\mU_3$ if $x_j$ is 0 or 1, respectively), or (ii) not restricted by the values of first-stage decisions $\bx \in \mX$ (since the constraints associated with $y_i$ become redundant when $\mU_2 = \emptyset$, and when $x_j = 1$ for all $(i,j) \in \mU_1$ and $x_j = 0$ for all $(i,j) \in \mU_3$). Practically, this implies that individual recourse variables not impacted by the first-stage decisions can be used to ``adapt" to new information. 

\subsection{Models with selective adaptability}\label{subsec:selective_adapt}

We give a few examples of ARBO problems with selective adaptability and then show that our examination of this problem structure comes without loss of generality. We begin with two examples of sequential decision-making problems that naturally exhibit this linking constraint structure:

\begin{itemize}
    \item In a facility location problem, a planner can only make an assignment (or shipping) decision between a potential location $j$ and an individual customer $l$ if a facility (e.g., hub, warehouse, factory) has been placed at location $j$ \citep{kammerling2020oracle}. Assuming that customer demand is given by $\bxi$, the sequential nature of this decision-making process can be represented by linking constraints of the form $y_{lj} \leq x_j$, where $y_{lj}$ is a particular location-customer assignment decision and $x_j$ determines whether a facility is placed in location $j$. 

    \item In a project investment problem, an investor with a limited budget aims to maximize investment returns across two decision stages by investing in different projects where early-stage investments come with risk but higher return rates \citep{arslan2022decomposition}. For each project $i \in \{1,\ldots,n\}$, the linking constraint $y_i \geq x_i$ captures early-stage commitment decisions. We revisit this problem in our numerical experiments in Section \ref{sec:experiments}.

\end{itemize}

The concept of selective adaptability can also be used to as a framework for introducing limited degrees of flexibility into static robust binary optimization problems. These models can be used to narrow down the set of choices to consider before model parameters are known exactly. Specifically, given the static robust problem
\begin{align*}
    \underset{\by \in \mY}{\text{min}} \ \, \underset{\bxi \in \Xi}{\text{max}} \quad \bxi^\top \by,
\end{align*}
we can define extensions where we have a \emph{budget} on the level of adaptability. For example, consider problem \eqref{model:1} with $\mS(\bx) = \{\by \in \mathbb{R}^{n=m} \, | \, \by \leq \bx\}$, $\bc = \bzero$, and $\mX = \{\bx \in \{0,1\}^{m} \, | \, \bw^\top \bx \leq \beta \}$ where $\bw$ and $\beta$ are non-negative. In this setting, $\beta$ controls the budget on adaptability, where higher values imply that more recourse decisions are ``available" in the second stage. For example, consider a daily route generation problem in a transportation context. If $\bw = \bone$ and $\beta = 0.1n$, where $n$ represents the number of links in a road network, the corresponding model would identify the most important 10\% of links (denoted by $\bx$) for constructing real-time vehicle routes (denoted by $\by \in \mY \cap \mS(\bx)$). Highlighting such subnetworks may be important for regular planning, communication, or training purposes. A similar problem is examined in the numerical experiments in Section \ref{sec:experiments}.

Finally, we emphasize that any ARBO problem without selective adaptability can be reformulated into one that has this property by introducing additional auxiliary variables into the recourse problem. Specifically, let $\bhF\bx + \bhG\by \leq \bhh$ be a set of linking constraints that are not in the form of equation \eqref{cons:linking} and let $\hat\mI \subseteq \{1, \ldots, m\}$ be the indices of the $\bx$ variables that appear in these constraints. Now, let $\boldsymbol{y^\text{aux}} \in \{0,1\}^{|\hat\mI|}$ be a new set of auxiliary recourse variables. For each $j \in \hat\mI$, we can replace $x_j$ in $\bhF\bx + \bhG\by \leq \bhh$ with the corresponding $y^\text{aux}_j$ variable, making these constraints a part of $\mY$,  and redefine the linking constraints as the set $\mS(\bx) = \{y^\text{aux}_j = x_j, \, \forall j \in \hat\mI\}$. This reformulation technique is inspired by \cite{arslan2022decomposition}, and implies that all results derived in our paper are applicable to any ARBO problem, with or without selective adaptability.

\subsection{Preliminary model insights}\label{subsec:prelim_insights}

Problem \eqref{model:1}, which is a min-max-min optimization problem, can be rewritten as the infinite-dimensional MILP model
\begin{subequations}\label{model:mono1}
\begin{alignat}{2}
 \underset{v, \bx, \by^\xi}{\text{min}} \quad & \bc^\top \bx + v \\
\text{s.t.} \quad 
& \bxi^\top \by^\xi \leq v, \quad && \forall \bxi \in \Xi \label{mono1:cons1}\\
& \by^\xi \in \mY , \quad && \forall \bxi \in \Xi \label{mono1:cons2}\\
& \by^\xi \in \mS(\bx), \quad && \forall \bxi \in \Xi \label{mono1:cons3}\\
& \bx \in \mX,
\end{alignat}
\end{subequations}
where $v$ represents the worst-case recourse objective value, and where a set of variables ($\by^\xi$) and constraints (\eqref{mono1:cons1} - \eqref{mono1:cons3}) must be defined for every $\bxi \in \Xi$. 

\cite{zeng2013solving} proposed a nested constraint-and-column generation to iteratively refine an approximation of this infinite-dimensional MILP formulation. In their approach, the set $\Xi$ in model \eqref{model:mono1} is replaced with a finite set $\hat\Xi$ which is iteratively enlarged by solving a subproblem that computes points in $\Xi$ to add to the set $\hat\Xi$. However, since model \eqref{model:mono1} is a complex combinatorial problem even for very small sets $\hat\Xi \subset \Xi$, the approach, which adds an additional copy of variables and constraints to the model at each iteration, can quickly become intractable. Furthermore, as discussed in Section \ref{subsec:RO_litreview}, the nested constraint-and-column generation algorithm may be challenging to implement, and finite convergence may be neither quick nor guaranteed.

In this paper, we begin by showing that model \eqref{model:mono1} permits a more tractable MILP reformulation when $\mS(\bx)$ satisfies selective adaptability. This reformulation is shown in the next proposition, which relies on the following lemma. All proofs can be found in the Electronic Companion. 

\begin{lem}\label{propn:conv_seperation}

$\conv(\mY \cap \mS(\bx)) = \conv(\mY) \cap \mS(\bx), \ \  \forall \bx \in \{0,1\}^m$. 

\end{lem}

\begin{propn}\label{thm1}
If Problem \eqref{model:1} has selective adaptability, then model \eqref{model:mono1} can be reformulated as shown, where $\by$ becomes a vector of continuous variables: 
\begin{subequations}\label{model:mono_simp}
\begin{align}
\underset{v, \bx, \by}{\mathrm{min}} \quad & \bc^\top\bx + v \\
\mathrm{s.t.} \quad & \bxi^\top \by \leq v, \quad \forall \bxi \in \Xi \label{model:mono_simp_cons1}\\
& \by \in \conv(\mY) \label{model:mono_simp_cons2}\\
& \by \in \mS(\bx) \label{model:mono_simp_cons}\\
& \bx \in \mX.
\end{align}
\end{subequations}
\end{propn}

Proposition \ref{thm1} is important for two main reasons, and motivates the rest of this paper. First, model \eqref{model:mono_simp} is equivalent to a static RO problem, i.e., where both $\bx \in \{0,1\}^m$ and $\by \in \mathbb{R}^n$ represent ``first-stage" decision variables. Practically, this means that we can solve the original ARBO problem as a MILP model if we have a polyhedral representation of $\conv(\mY)$. Section \ref{sec:BDD} focuses on the use of network flow constraints in an extended space to describe $\conv(\mY)$. Second, Proposition \ref{thm1} provides intuition for developing approximation schemes. In particular, model \eqref{model:mono_simp} implies that replacing $\conv(\mY)$ with an inner-(or outer-) approximation will generate a valid primal (or dual) bound on the optimal value of \eqref{model:mono_simp}. For example, replacing constraint \eqref{model:mono_simp_cons2} with the continuous relaxation of $\mY$, defined as $\rel(\mY) := \{\by \in [0,1]^n \, | \, \bG\by \geq \bh\}$, generates a lower bound on the optimal value of model \eqref{model:mono_simp}. Approximation schemes will be the focus of Section \ref{sec:approximations}.

Finally, there may naturally exist problems where the relaxation of $\mY$ is an integral polytope, i.e., $\rel(\mY) = \conv(\mY)$. One example is the budgeted adaptive routing problem mentioned in Section \ref{subsec:selective_adapt}, if we assume that the routing problem is a shortest path problem. In this setting, the adaptive robust problem can solved directly as a monolithic MILP model. 

\begin{propn}
If $\rel(\mY) = \conv(\mY)$, then model \eqref{model:mono1} is equivalent to
\begin{equation}
\begin{aligned}\label{model:compact_MIP}
\underset{\bx, \by, \blambda}{\mathrm{min}} \quad & \bc^\top \bx + \bd^\top \blambda \\
\mathrm{s.t.} \quad & \bT^\top \blambda = \by\\
& \bG \by \geq \bh\\
& \by \in \mS(\bx)\\
& \bx \in \mX, \ \ \by \in [0,1]^n, \ \ \blambda \geq \mathbf{0}
\end{aligned}
\end{equation}

by reformulating the uncertainty set using duality conditions \citep{gorissen2015practical}.
\end{propn}

We conclude this section with two remarks on the generalizability of the results in this paper when certain modeling assumptions are relaxed.

\begin{rem}
First-stage decisions can be mixed-integer, i.e., $\mX \subseteq \{0,1\}^{m_1} \times \mathbb{Z}^{m_2} \times \mathbb{R}^{m_3}$, as long as only the binary variables defining $\mX$ appear in the linking constraints $\mS(\bx)$. 
\end{rem}

\begin{rem}
The uncertainty set $\Xi$ can be any convex set. For example, if $\Xi$ is an ellipsoidal uncertainty set, then the MILP formulations we derive will become mixed-integer quadratic programs. 
\end{rem}

\section{Constrained Network Flow Reformulations}\label{sec:BDD}

In this section, we present a MILP formulation for model \eqref{model:mono_simp} when $\rel(\mY) \neq \conv(\mY)$. Specifically, in Section \ref{subsec:BDD_reform} we introduce decision diagrams and the corresponding network flow formulation used to obtain a polyhedral description of $\conv(\mY)$. Then, in Section \ref{subsec:BDD_final}, we integrate the formulation with first-stage decisions.

\subsection{Reformulating the recourse feasible space}\label{subsec:BDD_reform}

We first describe a general procedure for obtaining a decision diagram encoding of the feasible solutions in $\mY$, then represent this diagram using a network flow formulation.

\subsubsection{Decision diagram formulation.}\label{subsubsec:DDformulation}

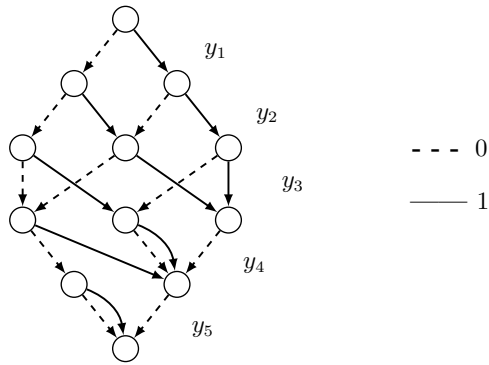
\begin{figure}[h]
\centering
\resizebox{0.4\textwidth}{!}{
 \begin{tikzpicture}[-latex]
 
   \node[draw, semithick,circle, minimum size = 0.4cm] (s1) at (0,0) {};
   \node[draw, semithick,circle, minimum size = 0.4cm] (s2n1) [below left = 0.7cm and 0.5cm of s1] {};
   \node[draw, semithick, circle, minimum size = 0.4cm] (s2n2)  [below right = 0.7cm and 0.5cm of s1] {};
   \node[draw, semithick,circle, minimum size = 0.4cm] (s3n2) [below right = 0.7cm and 0.5cm of s2n1] {};
   \node[draw, semithick, circle, minimum size = 0.4cm] (s3n1) [below left = 0.7cm and 0.5cm of s2n1] {};
    \node[draw, semithick, circle, minimum size = 0.4cm] (s3n3) [below right = 0.7cm and 0.5cm of s2n2] {};
    \node[draw, semithick, circle, minimum size = 0.4cm] (s4n1) [below = 0.7cm of s3n1] {};
    \node[draw, semithick, circle, minimum size = 0.4cm] (s4n2) [below = 0.7cm of s3n2] {};
    \node[draw, semithick, circle, minimum size = 0.4cm] (s4n3) [below = 0.7cm of s3n3] {};
    \node[draw, semithick, circle, minimum size = 0.4cm] (s5n1) [below right = 0.7cm and 0.5cm of s4n1] {};
    \node[draw, semithick, circle, minimum size = 0.4cm] (s5n2) [below right = 0.7cm and 0.5cm of s4n2] {};
    \node[draw, semithick, circle, minimum size = 0.4cm] (s6n1) [below right = 0.7cm and 0.5cm of s5n1] {};
    

    \node [right = 2.5cm of s3n3] (0)  {\textbf{- - -} $\, \, 0$};
    \node [below = 0.3cm of 0] (1)  {\textbf{------} $1$};

    \path[draw, line width = 0.3mm]
    (s1) edge node [xshift = 1cm] {$ y_1 $} (s2n2)
    (s2n1) edge node [yshift = -0.3cm] {$ $} (s3n2)
    (s2n2) edge node [xshift = 1cm] {$y_2$} (s3n3)
    (s3n1) edge node [xshift = 0.3 cm] {$ $} (s4n2)
    (s3n2) edge node [xshift = 0.3cm] {$ $} (s4n3)
    (s3n3) edge node [xshift = 1cm] {$y_3$} (s4n3)
    (s4n1) edge node [xshift = 0.3cm] {$ $} (s5n2)
    (s4n2) edge [bend left] node [xshift = 1cm] {$ $} (s5n2)
    (s5n1) edge [bend left] node [xshift = 1cm] {$ $} (s6n1);
    
    \path[draw, dashed, line width = 0.3mm]
    (s1) edge node [xshift = 0.3cm] {$ $} (s2n1)
    (s2n1) edge node [xshift = 0.3cm] {$ $} (s3n1)
    (s2n2) edge node [xshift = 0.3cm] {$ $} (s3n2)
    (s3n1) edge node [xshift = 0.3cm] {$ $} (s4n1)
    (s3n2) edge node [xshift = 0.3cm] {$ $} (s4n1)
    (s3n3) edge node [xshift = 0.3cm] {$ $} (s4n2)
    (s4n1) edge node [xshift = 0.3cm] {$ $} (s5n1)
    (s4n2) edge node [xshift = 0.3cm] {$ $}  (s5n2)
    (s4n3) edge node [xshift = 0.8cm, yshift = -0.2cm] {$y_4$} (s5n2)
    (s5n1) edge node [xshift = 0.3cm] {$ $} (s6n1)
    (s5n2) edge node [xshift = 0.8cm, yshift = -0.2cm] {$y_5$} (s6n1);

\end{tikzpicture}}
\vspace{0.5cm}
\caption{A decision diagram with six layers, where the zero-one label of each arc is indicated using a dashed or solid line.}\label{fig:BDD_ex}
\end{figure}

A binary decision diagram (BDD) is a graphical structure that can be used to encode the feasible space of a binary optimization problem \citep{bergman2016decision}.  Specifically, a BDD is a directed acyclic graph $\mD = (\mV, \mA)$ with nodes $\mV$ and arcs $\mA$. The nodes are partitioned into $n+1$ non-empty layers $\mV = (\mV_1, \ldots, \mV_{n+1})$, while the directed arcs connect nodes in adjacent layers from $\mV_i$ to $\mV_{i+1}$ for $i \in \{1, \ldots, n\}$. The sets $\mV_1 = \{\texttt{r}\}$ and $\mV_{n+1} = \{\texttt{t}\}$ are each composed of a single node, defined as the root node and terminal node, respectively. Each arc in the network has a label of either zero or one, and $\mA^0_i$ and $\mA^1_i$ define the subset of zero and one arcs leaving nodes in $\mV_i$, respectively. A decision diagram $\mD$ is a valid representation of a feasible space $\mY$ if each path from root node \texttt{r} to terminal node \texttt{t} in $\mD$ can be mapped to a solution $\by \in \mY$, and vice versa. This mapping of path to solution is defined by the zero-one label on each arc in the path. Specifically, for any arc $j$ in the path, if $j \in \mA^0_i$ then $y_i = 0$, or if $j \in \mA^1_i$ then $y_i = 1$. For example, a decision diagram of $\mY := \{\by \in \{0,1\}^5  \, | \, y_1 + y_2 + 2y_3 + 2y_4 + 3y_5 \leq 4 \}$ is given in Figure \ref{fig:BDD_ex}. Finally, we note that the definition of BDDs in the literature generally includes arc weights that correspond to the value of objective coefficients of an optimization problem. Since these coefficients are not deterministic in our problem, we instead consider only ``unweighted" decision diagrams.

Decision diagrams are closely related to the state-transition graph in the dynamic programming literature \citep{hooker2013decision}. In particular, nodes and arcs in the decision diagram can be mapped to ``states" and feasible ``actions" of a recursive formulation where decisions are made sequentially. Many binary optimization problem structures admit simple recursive formulations that can be used to obtain decision diagrams; we refer to \cite{bergman2016decision}, \cite{bergman2022network} and \cite{de2022arc} for a comprehensive summary of recursive formulations for a variety of problem structures. In general, a recursive formulation of a deterministic binary optimization problem can be written as a Bellman equation of the form
\begin{align}
V_i(\bS) = \underset{y_i \in Q_i(\bS)}{\text{max}} \ \  \bigg \{f_i(\bS, y_i) + V_{i+1}(T_i(\bS,y_i)) \bigg \},
\end{align}
where $\bS$ denotes the state of the system, $\mQ_i(\bS) \subseteq \{0,1\}$ denotes the set of feasible actions at stage $i$ for state $\bS$, $T_i(\bS, y_i)$ defines the new state after taking action $y_i$, and the pair $V_i(\bS)$ and $f_i(\bS, y_i)$ are used to capture the long-term values and immediate rewards of taking specific actions, respectively. However, since we only need to generate unweighted decision diagrams, we only need the state-transition graph corresponding to the states $\bS$ and feasible actions $Q_i(\bS)$. To define the diagram, we first modify $T_n(\bS, y_i)$ to be $T_n(\bS, y_i) = \{\texttt{t}\}$, and then map each state to a node and each feasible action $y_i \in Q_i(\bS)$ to a unique arc with label $y_i$ that links $\bS$ and $T_i(\bS, y_i)$. This creates a diagram with a single root node \texttt{r} and terminal node \texttt{t}, where each path from \texttt{r} to \texttt{t} corresponds to a sequence of actions in the recursive formulation.

\begin{ex} \label{ex:BDD_construction}
Consider a feasible set $\mY = \big\{ \by \in \{0,1\}^n \, \big | \, \sum_{i =1}^n g_iy_i \leq h \big \}$ that is defined by a knapsack constraint. The recursive formulation of $\mY$ is defined by feasible actions $\mQ_i(\bS) = \{y_i \in \{0,1\} \, | \, \bS + g_iy_i \leq h\}$, state-transition function $T_i(\bS, y_i) = \bS + g_iy_i$ with an initial state $\bS_1 = 0$. Each stage $i$ in the recursive formulation thus corresponds to a layer $i$ in the diagram. Each state in layer $i$ represents the total amount of capacity used by the selection of items among $1,\ldots, i$ and feasible actions correspond to whether an additional item $i+1$ can be placed in the knapsack given the current state. To obtain a single terminal node in the decision diagram, we merge all nodes in layer $n+1$, which is equivalent to replacing $T_n(\bS, y_n) = \bS + g_ny_n$ with $\mT_n(\bS, y_n) = h$ in the recursive formulation. Note that this change in the state-transition matrix at the final stage $n$ does not change the set of feasible actions in the recursive formulation.
\end{ex}

Finally, we note for any given decision diagram, it may be possible to obtain a \emph{reduced} version that encodes the same set of feasible solutions in a much fewer number of nodes and arcs. Figure \ref{fig:BDD_ex} is a reduced decision diagram for the knapsack set defined earlier. To obtain reduced decision diagrams, we can use a simple bottom-up merging technique \citep{bryant1992symbolic, bergman2016decision}. Starting with nodes in layer $n$, we can merge any nodes which have the same set of outgoing arc types and destinations, repeating this procedure from layers $n-1$ to $1$.

\subsubsection{Recourse network flow formulation.}

Let $\mD$ denote a decision diagram representation of $\mY$. Following the notation of \cite{castro2022decision}, we use $\text{NF}(\mD)$ to denote the network flow model of $\mD$, which relates arc flows $\bz$ to values of $\by \in \mY$ in the recourse problem. Specifically, let $\NF(\mD)$ be defined as
\begin{align}
\NF(\mD) := \bigg \{(\by, \bz) \in \mathbb{R}^n \times \mathbb{R}^{|\mA|}_+ \, \bigg | \, \bA \bz = \bb, \ y_i = \sum_{j \in \mA^1_i} z_j, \ \forall i \in \{1, \ldots, n\} \bigg\},
\end{align}
where $\bA \in \{-1,0,1\}^{|\mV| \times |\mA|}$ denotes the node-arc incidence matrix and $\bb$ a vector of zeros with the exception of $b_1 = -1$ and $b_{|\mV|} = 1$. The constraint set $\bA\bz = \bb$ defines standard flow conservation constraints. The second set of constraints link the value of each variable $y_i$ to the sum of total flow over the one-arcs $\mA_i^1$ in layer $i \in \{1,\ldots, n\}$.

\subsection{A complete network flow reformulation}\label{subsec:BDD_final}

A key property of $\NF(\mD)$ is that its projection onto the variables $\by$, denoted by $\proj_\by(\NF(\mD))$, is equal to the convex hull of $\mY$ \citep{castro2022decision}. We can use this property to derive an MILP formulation of Problem \eqref{model:1}, as shown in the next proposition.

\begin{propn}
Problem \eqref{model:1} can be reformulated into the constrained network flow problem
\begin{subequations}\label{model:mono_adapt}
\begin{align}
\underset{\bx, \by, \bz, \blambda}{\mathrm{min}} \quad & \bc^\top \bx + \bd^\top \blambda \\
\mathrm{s.t.} \quad & \bT^\top \blambda = \by\\
& y_i = \sum_{j \in \mA^1_i} z_j \quad \forall i \in \{1, \ldots, n\}\\
& \bA \bz = \bb \\
& \by \in \mS(\bx) \\
& \bx \in \mX, \ \bz \geq \bzero.
\end{align} 
\end{subequations}
\end{propn}

Model \eqref{model:mono_adapt} is thus an exact MILP reformulation of any ARBO problem with selective adaptability. As we will show in Section \ref{sec:experiments}, this model can be tractably solved for problems of smaller sizes.

\begin{ex}\label{ex:2stage_ex1}
Consider the adaptive robust knapsack problem
\begin{align}\label{model:BDD_ex}
\begin{split}
\underset{\bx \in \{0,1\}^5}{\textup{min}} \ \  \underset{\bxi \in \Xi}{\textup{max}} \ \ \underset{\by}{\textup{min}} \ \quad & \bc^\top \bx + \bxi^\top \by \\
& \by \leq \bx\\
& \by \in \mY := \left \{\by \in \{0,1\}^5 \, \Big | \, y_1 + y_2 + 2y_3 + 2y_4 + 3y_5 \leq 4 \right \}.
\end{split}
\end{align}
Recall that the decision diagram for $\mY$ is illustrated in Figure \ref{fig:BDD_ex}. Suppose $\Xi = \left \{ \bxi \,  | \, \bxi \geq \bxi^0, \, \norm{\bxi - \bxi^0}_1 \leq \delta \right \}$. Then, problem \eqref{model:BDD_ex} can be reformulated into model \eqref{model:mono_adapt}, where $\bA$ is the node-arc incidency matrix of the decision diagram given in Figure \ref{fig:BDD_ex}, and 
\begin{align*}
\bT = \begin{bmatrix} \ 
-1 & 0 & 0 & 0 & 0\\
0 & -1 & 0 & 0 & 0\\
0 & 0 & -1 & 0 & 0\\
0 & 0 & 0 & -1 & 0\\
0 & 0 & 0 & 0 & -1\\
1 & 1 & 1 & 1 & 1 \ 
\end{bmatrix}, \ \ 
\bb = \begin{bmatrix} \ 
-1\\
0\\
0\\
\vdots\\
0\\
1 \ 
\end{bmatrix}, \ \ 
\bd = \begin{bmatrix} \ 
-\xi^0_1\\
-\xi^0_2\\
-\xi^0_3\\
-\xi^0_4\\
-\xi^0_5\\
(\bxi^0)^\top \mathbf{1} + \delta \ 
\end{bmatrix}, \quad 
\begin{matrix*}[l]  
y_1 = z_2 & \\
y_2 = z_4 + z_6 & \\
y_3 = z_8 + z_{10} + z_{12} & \\
y_4 = z_{14} + z_{16} & \\
y_5 = z_{19}. & 
\end{matrix*}
\end{align*}

The indices $\{2, 4, 6, 8, {10}, {12}, {14}, {16}, {19}\}$ correspond to the solid lines in Figure \ref{fig:BDD_ex} when labeled from left to right in each layer, starting with the first layer. 
\end{ex}

\section{Network Flow Approximations}\label{sec:approximations}

In this section, we use approximate decision diagrams to derive new approximation methods for large-scale ARBO problems. We first define restricted and relaxed decision diagrams, which we then use to propose compact network flow models that generate first-stage solutions, primal bounds and dual bounds, respectively. We then present a multi-network flow model, which serves as a more general framework for generating dual bounds. Finally, we conclude the section by outlining a procedure for evaluating the quality of any computed first-stage solution.

\subsection{Approximate decision diagrams}

A restricted decision diagram of $\mY$ contains paths that map to a strict subset of the feasible solutions in $\mY$. On the other hand, a relaxed decision diagram contains paths that map to a superset of solutions that include all solutions in $\mY$. Restricted and relaxed decision diagrams can be generated by merging states in the recursive formulation of $\mY$ \citep{castro2022decision}. Next, we outline two general ``top-down" approaches for building restricted and relaxed decision diagrams.

\subsubsection{Merging with width-based thresholds.}\label{subsubsec:width-based}

In Algorithm \ref{algo:RC} we present a common approach that merges nodes whenever a threshold on the diagram's ``width" has been exceeded. Specifically, each layer in the diagram is constrained to have a width of at most $W$, that is, there can be at most $W$ nodes per layer. Considering the recursive formulation, this width constraint is equivalent to having at most $W$ number of different states at stage $i$. Note that if $W = \infty$ in Algorithm \ref{algo:RC}, then an exact decision diagram will be generated. 

\begin{algorithm}[h]
\caption{\footnotesize Width-Based Decision Diagram Construction Procedure}
\label{algo:RC}
\vspace{0.1cm}
\footnotesize \textbf{Input:} A recursive formulation of $\mY$, a width parameter $W$\\
\textbf{Output:} A decision diagram $\mD$
\begin{algorithmic}[1]
\State Initialization: let $\bS_1(\texttt{r})$ denote the initial state of the system (at root node $\texttt{r}$), \, $\mV_1 = \{\texttt{r}\}, \newline \mV_2, \ldots, \mV_n = \emptyset, \, \mV_{n+1} = \{\texttt{t}\}, \, \mA = \emptyset$
\For{$i \in \{1, \ldots, n-1\}$}
\For{$u \in \mV_i$} \label{alg_step:init}
\For{$y_i \in Q_i(\bS(u))$} 
\If{$ \exists \ u' \in \mV_{i+1}$ where $\bS(u') = T_i(\bS, y_i)$} add arc from $u$ to $u'$ with label $y_i$
\Else{ \ add node $u'$ with state $T_i(\bS, y_i)$ to $\mV_{i+1}$ and add arc from $u$ to $u'$ with label $y_i$} \label{alg_step:next_layer}
\EndIf
\EndFor
\EndFor
\If{$|\mV_{i+1}| > W$} 
\State $\hat\mV_{i+1} \gets \textsf{select}(\mV_{i+1}, |\mV_{i+1}| - W + 1)$\label{alg_step:select}
\State $\textsf{merge}(\hat\mV_{i+1})$ \label{alg_step:merge_discard}
\EndIf
\EndFor
\For{$u \in \mV_n$} \label{alg_step:final_layer}
\For{$y_i \in Q_i(\bS(u))$}
add arc from $u$ to $\texttt{t}$ with label $y_i$ 
\EndFor
\EndFor
\State \textbf{Reduce} $\mD = (\mV, \mA)$ \label{alg_step:reduce}\\
\Return $\mD$ \label{alg_step:return} 
\end{algorithmic}
\end{algorithm}

In step \ref{alg_step:select}, the algorithm requires a rule for selecting $(|\mV_{i+1}| - W + 1)$ nodes from the set $\mV_{i+1}$. Random selection is the simplest and most commonly used procedure. Nonetheless, designing node selection procedures have become a more active topic of research in recent years \citep[e.g.,][]{van2022graph}, although such procedures focus on specific classes of deterministic problems. 

Once this subset of nodes $\hat\mV_i$ has been selected, they can be merged to create either restricted or relaxed decision diagrams. To create relaxed decision diagrams, merging operators must assign a new state to each merged node such that all subsequent feasible actions in the original recursive formulation remain feasible under this new state \citep{hooker2013decision}. 
Similarly, a restricted decision diagram can be created by assigning a new state to the merged node such that a subset of feasible actions in the original recursive are retained, without introducing infeasible actions. We give an example below.

\begin{ex}\label{ex:knapsack_constraint}
Consider a feasible space defined by a single knapsack constraint (see Example \ref{ex:BDD_construction}). A restricted decision diagram can be generated by merging the nodes $\hat\mV_{i+1}$ and assigning the new node a state that is the minimum value of the states in the merged set. A relaxed decision diagram can be generated by assigning the new node the maximum value of the states in the merged set.
\end{ex}

We note that the node merge operation can be performed while building a layer, rather than after the layer is completely built, to reduce the memory requirement, if desired. Furthermore, we note that we can also generate restricted diagrams simply by discarding nodes rather than merging them in Algorithm \ref{algo:RC}, since we are effectively removing all feasible actions associated with that state-stage combination in the recursive formulation. Finally, step \ref{alg_step:reduce} is to reduce the size of the decision diagram using the bottom-up approach described at the end of Section \ref{subsubsec:DDformulation}.


\subsubsection{Merging with distance-based thresholds.}\label{subsubsec:distance-based} We now propose a merging approach based on the similarity of states in each layer, that is, we merge two nodes only when their respective state values are within some ``distance" of each another. Algorithm \ref{alg:QDD} outlines this approach, which revolves around a partitioning of nodes in each layer into subgroups, where nodes in each subgroup are then merged. Any suitable set partitioning method can be used, as long two conditions are met: (i) within each subgroup, some pre-defined notion of distance between any pair of states does not exceed a user-specified parameter $Q$, and (ii) for any two subgroups, there exists a pair of nodes, one node in each subgroup, where the distance exceeds $Q$. The latter condition ensures the fewest number of partitions given $Q$. To the best of our knowledge, we are the first to outline and implement this distance-based approach for generating approximate decision diagrams. 

\begin{algorithm}[h]
\caption{\footnotesize Distance-Based Decision Diagram Construction Procedure}\label{alg:QDD}
\vspace{0.1cm}
\footnotesize \textbf{Input:} A recursive formulation of $\mY$, a distance function $d(\cdot,\cdot)$ and distance parameter $Q$\\
\textbf{Output:} A decision diagram $\mD$
\begin{algorithmic}[1]
\State Initialization: let $\bS_1(\texttt{r})$ denote the initial state of the system (at root node $\texttt{r}$), \, $\mV_1 = \{\texttt{r}\}, \newline \mV_2, \ldots, \mV_n = \emptyset, \, \mV_{n+1} = \{\texttt{t}\}, \, \mA = \emptyset$
\For{$i \in \{1, \ldots, n-1\}$}
\State Steps \ref{alg_step:init} to \ref{alg_step:next_layer} from Algorithm \ref{algo:RC}
\State  $\hat\mV_{i+1, 1}, \hat\mV_{i+1, 2}, \ldots \gets \textsf{partition}(\mV_{i+1}, Q)$ \label{alg_step:partition} 
\For{$\hat\mV_{i+1} \in \{\hat\mV_{i+1, 1}, \hat\mV_{i+1, 2}, \ldots \}$}
\State $\textsf{merge}(\hat\mV_{i+1})$
\EndFor
\EndFor
\State Steps \ref{alg_step:final_layer} to \ref{alg_step:return} from Algorithm \ref{algo:RC}
\end{algorithmic}
\end{algorithm}

\begin{ex}\label{ex:partition_function}
We give an example of a \textsf{partition} function (step \ref{alg_step:partition} of Algorithm \ref{alg:QDD}). Consider the feasible space defined in Example \ref{ex:knapsack_constraint}, for which we can define the following node partitioning procedure: (i) order nodes according to state values, from smallest to largest, (ii) starting with the smallest state, create a new subgroup for the next node if and only if its state is more than $Q$ units larger than the smallest state in the current subgroup. Here, distance between states is simply defined as the absolutely difference between state values. 
\end{ex}

The distance-based merging approach is motivated by our initial attempts at using width-based merging in preliminary experiments. We found that it was difficult to finely tune the size of the diagram and quality of the approximation through the use of the width threshold parameter $W$. In particular, when the width threshold is exceeded, it is often exceeded by a large margin, which typically results in the merging of hundreds or thousands of nodes into a single node. Furthermore, because the width constraint must be satisfied at each layer, there are essentially no restrictions on which nodes can or cannot be merged. In contrast, the key advantage of the distance-based approach is that both the size and the quality of the approximation is closely tied to the value of $Q$. In distance-based merging, $Q$ determines precisely which nodes in each layer can and cannot be merged based on the similarity of their states. Furthermore, at each layer, there may be many subgroups but the number of nodes to be merged within each subgroup may be small. In general, we find that this distance-based merging procedure leads to approximations of higher quality and more precision in tuning the size and tractability of the corresponding models.

Finally, we remark that the distance-based approach can be further customized to generate even more refined approximations, for example, by merging subgroups with some probability $p < 1$. This additional feature could help generate approximations that are in-between those generated solely by incrementally increasing the value of $Q$ (which could be discrete, for example, when all states are integer-valued as in Example \ref{ex:knapsack_constraint}).

\subsection{Primal bounds}\label{subsec:primal_bound}

For a restricted diagram $\mD_{\texttt{inner}}$, it is the case that $\proj_\by(\text{NF}(\mD_{\texttt{inner}})) \subseteq \conv(\mY)$ since $\mD_{\texttt{inner}}$ contains a subset of the feasible solutions of $\mY$. Thus, restricted decision diagrams can be used to design network flow models that generate primal bounds on the ARBO problem. 

\begin{propn}
Let $\mD_{\texttt{inner}}$ denote a restricted decision diagram of $\mY$. Then, the model 
\begin{subequations} \label{model:inner_approx}
\begin{align} 
\underset{\bx, \by, \blambda}{\mathrm{min}} \quad & \bc^\top \bx + \bd^\top \blambda \\
\mathrm{s.t.} \quad & \bT^\top \blambda = \by\\
& \by \in \proj_\by(\NF(\mD_{\texttt{inner}})) \\
& \by \in \mS(\bx) \\
& \bx \in \mX
\end{align} 
\end{subequations}
generates a feasible solution $\bx \in \mX$ and an upper bound on the optimal objective value of model \eqref{model:1}. 
\end{propn}

There are conceptual connections between model \eqref{model:inner_approx} and $K$-adaptability in that both approximation methods are derived by restricting the set of possible recourse decisions. Nonetheless, there are two important differences. First, the approximation scheme defined by model \eqref{model:inner_approx} relies on continuous recourse variables, while $K$-adaptability relies on discrete recourse variables. Second, the size of model \eqref{model:inner_approx} can be more precisely controlled by incrementally changing the width of the decision diagram. On the other hand, in $K$-adaptability, each incremental increase in the value of $K$ requires adding a set of recourse decision variables $\by^{K+1}$ and corresponding constraints $\mY$. As we will illustrate in our numerical experiments (Section \ref{sec:experiments}), problems with even a few $K$ can quickly become intractable. 


\subsection{Dual bounds}\label{subsec:dual_bounds}

A relaxed decision diagram $\mD_{\texttt{outer}}$ of $\mY$ represents a superset of $\mY$, which implies that $\conv(\mY) \subseteq \proj_\by(\text{NF}(\mD_{\texttt{outer}}))$. Relaxed decision diagrams can thus be used to derive a dual bound. 

\begin{propn}
Let $\mD_{\texttt{outer}}$ denote a relaxed decision diagram of $\mY$. Then, the model
\begin{subequations}\label{model:outer_approx}
\begin{align} 
\underset{\bx, \by, \blambda}{\mathrm{min}} \quad & \bc^\top \bx + \bd^\top \blambda \\
\mathrm{s.t.} \quad & \bT^\top \blambda = \by\\
& \by \in \proj_\by(\NF(\mD_{\texttt{outer}})) \\
& \by \in \rel(\mY)\\
& \by \in \mS(\bx) \\
& \bx \in \mX
\end{align} 
\end{subequations}
generates a feasible solution $\bx \in \mX$ and a lower bound on the optimal objective value of model \eqref{model:1}.
\end{propn}

Model \eqref{model:outer_approx} provides a dual bound that is at least as strong as that with a simple continuous relaxation of recourse decisions $\by$. However, this bound can be stronger, since in general, $\rel(\mY) \nsubseteq \proj_\by(\NF(\mD_{\texttt{outer}}))$ and $\proj_\by(\NF(\mD_{\texttt{outer}})) \nsubseteq \rel(\mY)$. Specifically, $\rel(\mY)$ can have fractional extreme points, whereas $\proj_\by(\NF(\mD_{\texttt{outer}}))$ is an integral polyhedron but includes solutions $\bhy \in \{0,1\}^n$ that are not in $\mY$. The intersection of $\rel(\mY)$ and $\proj_\by(\NF(\mD_{\texttt{outer}}))$ can thus be a more accurate outer approximation of $\conv(\mY)$ than either of the two sets alone. 

While both models \eqref{model:inner_approx} and \eqref{model:outer_approx} are bounding techniques that rely on approximating $\conv(\mY)$ using a single decision diagram, we now propose a generalization where dual bounds can be derived with an outer approximation of $\conv(\mY)$ using a \emph{collection} of diagrams and feasible sets.

\begin{propn}\label{propn:constraint_decoupling}
Suppose that $\mY = \{\by \in \{0,1\}^n \, | \, (\bg^1)^\top \by \geq h_1, \ldots, (\bg^J)^\top \by \geq h_J\}$, and let $\mJ = \{1, \ldots, J\}$ denote the set of constraint indices. Now, suppose we are given a collection of subsets of indices $\mJ_1^1, \ldots, \mJ_k^1$ and $\mJ_1^2, \ldots, \mJ_q^2$ where $\cup_{i = 1}^k \mJ_i^1 \subseteq \mJ$ and $\cup_{i = 1}^q \mJ^2_i \subseteq \mJ$. For an arbitrary set $\mJ_i \subseteq \mJ$, let $\mD^{\mJ_i}$ and $\mD_\texttt{outer}^{\mJ_i}$ denote an exact and relaxed decision diagram for the feasible set described by the constraints with indices in $\mJ_i$, respectively. Then, the model 
\begin{equation}
\begin{aligned}
\label{model:multi-network}
\underset{\bx, \by, \blambda}{\mathrm{min}} \quad & \bc^\top \bx + \bd^\top \blambda \\
\mathrm{s.t.} \quad & \bT^\top \blambda = \by\\
& \by \in \proj_\by(\NF(\mD^{\mJ_i})), \quad && \forall \mJ_i \in \{\mJ_1^1, \ldots, \mJ_{k}^1\} \\
& \by \in \proj_\by(\NF(\mD^{\mJ_i}_{\texttt{outer}})), \quad && \forall \mJ_i \in \{\mJ_1^2, \ldots, \mJ_q^2\} \\
& \by \in \rel(\mY)\\
& \by \in \mS(\bx)\\
& \bx \in \mX
\end{aligned} 
\end{equation}
generates a feasible solution $\bx \in \mX$ and lower bound on the optimal objective value of model \eqref{model:1}.

\end{propn}

As a proof of concept, consider the model from Example \ref{ex:2stage_ex1}, but with an additional constraint of $y_1 + y_2 + y_3 \leq 2$. Rather than formulating an exact or relaxed decision diagram of
\begin{align*}
\mY := \Big\{ \by \in \{0,1\}^5 \, | \, y_1 + y_2 + y_3 \leq 2, y_1 + y_2 + 2y_3 + 2y_4 + 3y_5 \leq 4 \Big\},
\end{align*} 
we could generate two exact or relaxed decision diagram, one for each constraint, and combine them into the multi-network flow model \eqref{model:multi-network}.

In practice, multi-network approximations can be useful when single-network representations are too large or when the underlying problem structure does not admit an obvious recursive formulation. In these settings, we can generate exact or relaxed decision diagrams for subsets of constraints that do admit a simple recursive formulation. Furthermore, in many decision-making problems, a large subset of constraints defining $\mY$ may have a totally unimodular constraint coefficient matrix and integer right-hand sides, and thus define a feasible space for which its relaxation is an integral polytope. In this setting, we can generate decision diagrams only for the remaining constraints. The corresponding multi-network flow approximations may potentially be smaller in size, easier to implement, or better in quality than single-network approximations. We will further explore these models in the numerical experiments in Section \ref{subsec:assignment}.

\subsection{Evaluating the quality of a solution}\label{subsec:evaluation_oracle}

All approximation techniques presented in this section generate a feasible first-stage solution $\bx \in \mX$. We can compute the true objective value of this solution, which we denote using $z(\bhx)$, by solving
\begin{align}\label{ECmodel:evaluation_oracle}
z(\bhx) \, := \ \underset{\bxi \in \Xi}{\text{max}} \ \underset{\by \in \mY\cap \mS(\bhx)}{\text{min}} \quad & \bc^\top \bhx + \bxi^\top \by.
\end{align}

One way of solving for the value $z(\bhx)$ is to use the constraint generation method described in \cite{kammerling2020oracle}. Specifically, the authors consider an iterative algorithm between a master problem that computes a specific parameter realization $\hat{\bxi} \in \Xi$ and a subproblem which computes recourse solutions $\bhy \in \mY\cap \mS(\bhx)$. At iteration $k$, the master problem is 
\begin{align*}
\text{MP}(\bhy^{1},\ldots, \bhy^{k-1}) \ = \ \underset{v, \bxi}{\text{max}}  \ \ \Big \{ v \, \Big | \, v \leq \bc^\top \bhx + \bxi^\top \bhy^j \ \ \forall j \in \{1, \ldots, k-1\}, \ \bxi \in \Xi \Big \}.
\end{align*}
The optimal solution $\hat\bxi^k$ to this master problem is then passed to the subproblem
\begin{align*}
\text{SP}(\hat \bxi^k) \ = \ \underset{\by}{\text{min}} \ \ \Big\{ \bc^\top \bhx + (\hat\bxi^k)^\top \by \, \Big | \, \by \in \mY\cap \mS(\bhx) \Big\}.
\end{align*}
Let $\hat \by^{k}$ denote an optimal solution of the subproblem. If the objective value of the subproblem is less than the master problem objective value, then we add constraint $v \leq \bc^\top \bhx + \bxi^\top \hat\by^{k}$ into the master problem, and move to iteration $k+1$. Otherwise, $z(\bhx) = \bc^\top \bhx + (\hat{\bxi}^k)^\top \hat\by^k$ is the optimal objective value of Problem \eqref{ECmodel:evaluation_oracle}. Since the master problem is a linear program and the subproblem is one instantiation of the recourse problem, the effort required to compute the value of $z(\bhx)$ for a given $\bhx \in \mX$ is generally negligible relative to the general ARBO problem of solving for an optimal $\bhx \in \mX$.

Finally, recall that the approximation models presented in Section \ref{subsec:dual_bounds} generate dual bounds in addition to a solution $\bhx \in \mX$. For these models, we can thus compute a \emph{model-based optimality gap}, which we define as 
\begin{align}\label{eq:db_optgap}
    \frac{\text{dual bound} - z(\bhx)}{z(\bhx)}\cdot 100 \, ,
\end{align}
where `dual bound' denotes the optimal objective value of the approximation model. The model-based optimality gap serves as a proxy for the true optimality gap. In particular, the former is an upper bound on the latter. This distinction is important, because the true optimality gap can be challenging to compute in numerical experiments as it requires solving the ARBO problem exactly. Nonetheless, we will show in Section \ref{sec:experiments} that for most problems considered, the model-based optimality gap is in fact very low. This observation highlights that these approximation models not only compute high-quality solutions but can also verify the near-optimality of the solutions.

\section{Numerical Experiments}\label{sec:experiments}

In this section, we demonstrate and compare the effectiveness of three types of formulations: (i) exact network flow models, (ii) approximate single-network flow models, and (iii) approximate multi-network flow models. These models are applied to two ARBO problems, namely, a capital budgeting problem and a robust assignment problem. All experiments were coded in Python 3.7 and MILP models were solved using Gurobi 9.1.1 under default settings. Decision diagrams were created and manipulated using the NetworkX package. The experiments were conducted on an Macbook Pro (M1 Chip) with 16GB of RAM. Unless otherwise stated, all benchmark models that we discuss next were also implemented and solved under the same conditions.

\subsection{Capital budgeting problems}\label{subsec:capital_budgeting}

We first consider a capital budgeting problem, which seeks to compute a robust investment plan for a set of $n$ projects under an investment budget constraint. Several variants of this problem have been studied in the adaptive robust optimization literature \citep{hanasusanto2015k,subramanyam2020k, kammerling2020oracle, arslan2022decomposition, dumouchelle2023neur2ro}. In this problem, investment decisions are binary and can be made in two stages, i.e., projects can either be invested in at an early stage (denoted using $\bx \in \{0,1\}^n$) or in a later stage (denoted using $\by \in \{0,1\}^n$). The profitability of a project is not known with certainty in the first stage and is only revealed in the second stage. Early stage investors are rewarded with a first-mover advantage that entails a higher percentage on final profit generated by the project (otherwise, the optimal action would be to postpone all investment decisions until more information is revealed). 

For the rest of this section, we follow the model formulation from \cite{arslan2022decomposition} and use their publicly-available problem instances\footnote{see https://github.com/borisdetienne/RobustDecomposition (accessed February 2024)}. Specifically, we consider the following adaptive capital budgeting problem
\begin{subequations}
\begin{align}
\underset{\bx \in \{0,1\}^n}{\textup{max}} \ \  \underset{\bxi \in \Xi}{\textup{min}} \ \ \underset{\by}{\textup{max}} \ \quad & (1-f)(\bxi^\top \bx) + f(\bxi^\top \by) \\
\text{s.t.} \ \quad & \by \geq \bx \label{cons:CP_SA}\\
& \bg^\top \by \leq h \label{cons:CP_budget}\\ 
& \by \in \{0,1\}^n. \label{cons:CP_bounds}
\end{align}\label{model:CP}
\end{subequations}
In this model, $g_i$ denotes the cost of investing in project $i$, $h$ is the total investment budget, $f \in [0,1)$ captures the first-mover advantage, and $\xi_i$ is the payoff of project $i$ which is unknown in the first stage. The payoff of each project depends on a set of $M$ common risk factors $\alpha_1, \ldots, \alpha_m$. Specifically, it is assumed that $\xi_i = \sum_{j = 1}^M U_{ij}\alpha_j$ where $U_{ij} \in \mathbb{R}$ describes the impact that each risk factor $\alpha_j$ has on the payoff $\xi_i$. The risk factors $\alpha_j$ are assumed to reside in $[-1,1]$, i.e., 
\begin{align*}
\Xi = \left \{\bxi \, \middle |  \, \bxi = \bU\balpha, \balpha \in [-1,1]^M \right \}.
\end{align*}
Note that in this problem, constraint \eqref{cons:CP_SA} defines a selective adaptability constraint while constraints \eqref{cons:CP_budget}-\eqref{cons:CP_bounds} can be reformulated or approximated in an extended network space. 

We consider 300 instances of varying size where $n = \{10, 20, 30, 40, 50\}$ and where each value of $n$ is associated with 60 different instances \citep{arslan2022decomposition}. Each instance is characterized by a random sample of project costs and payoffs, as well as an investment budget defined as a fraction of total project costs, i.e., $h = m \sum_{i = 1}^n g_i$ with $m = \{0.2, 0.4, 0.6, 0.8\}$. 

Since \eqref{cons:CP_budget} is a knapsack constraint, we generate the exact network flow model by following the steps in Example \ref{ex:BDD_construction} and the reformulation technique outlined in Section \ref{sec:BDD}. We also examine approximate network flow models based on relaxed decision diagrams. These diagrams are generated using the distance-based merging approach discussed in Section \ref{subsubsec:distance-based}. Recall that in this approach, a user-specified value $Q$ bounds the distance between state values of nodes that are to be merged. In the capital budgeting problem, the state value of a node in layer $i$ defines the amount of weight that has been added to the knapsack based on decisions $y_1, \ldots, y_i$, and $Q$ represents the maximum difference in weight values between any two nodes to be merged. We follow the procedure described in Example \ref{ex:partition_function} to determine which nodes should be merged. Once merged, the new node is assigned the minimum state value of the nodes that were merged.

Since the approximation models generate first-stage solutions and dual bound, we can calculate optimality gaps for each solution. A key takeaway of the following numerical experiments is that these models can generate solutions that are verifiably near-optimal in very little time. 

\subsubsection{Model formulation and solution time.} We begin by examining the size, formulation time, and solution time of the exact and approximate network flow models. Table \ref{table:RDD_size} highlights the number of arcs in the decision diagrams used to represent constraint \eqref{cons:CP_budget}, shown over various values of $Q$. This number conveys the size of the corresponding network flow model, since each arc in the decision diagram corresponds to a continuous variable in our formulation. 

\begin{table}[h]
\centering
\setlength{\tabcolsep}{9pt}
{\renewcommand{\arraystretch}{1.2}%
\begin{tabular}{l|rrrrr}
 \toprule
 & $Q = 0$ & $Q = 1$ & $Q = 3$ & $Q = 5$ & $Q = 10$ \\
\midrule 
$n = 10$ & 166 & 164 (99\%) & 151 (91\%) & 134 (81\%) & 109 (65\%) \\
$n = 20$  & 2992 & 1972 (66\%) & 1198 (40\%) & 894 (30\%) & 554 (19\%) \\ 
$n = 30$  & 12334 & 6726 (55\%) & 3600 (29\%) & 2463 (20\%) & 1359 (11\%) \\ 
$n = 40$ & 28121 & 14621 (52\%) & 7374 (26\%) & 4912 (17\%) & 2622 (9\%) \\ 
$n = 50$  & 48285 & 24589 (51\%) & 12180 (25\%) & 7954 (16\%) & 4160 (9\%) \\
\bottomrule
\end{tabular}}
\caption{The average number of arcs in the reduced decision diagrams under different $Q$ values. The percentage of each value relative to that of the exact decision diagram (i.e., $Q = 0$) are given in parentheses.}\label{table:RDD_size}
\end{table}

Similarly, Table \ref{table:RDD_times} highlights the time it takes to generate the decision diagram and solve the network flow model. We make two observations. When $n$ increases, the time it takes to formulate and solve the exact model grows exponentially. This makes the exact models intractable beyond small values $n$. In contrast, as we increase $Q$ for a fixed value of $n$, the time it takes for both these processes can be reduced drastically. For example, when $n = 50$, the total time it takes to generate and solve the network flow models can be reduced by several orders of magnitude when $Q \geq 3$. 

\begin{table}[h]
\centering
\setlength{\tabcolsep}{18pt}
{\renewcommand{\arraystretch}{1.2}%
\begin{tabular}{c|rrrrr}
 \toprule
 & $Q = 0$ & $Q = 1$ & $Q = 3$ & $Q = 5$ & $Q = 10$ \\
\midrule 
\emph{Instances} & \multicolumn{5}{l}{\emph{Average build time for the reduced diagram, in seconds}} \\ 
 $n = 10$  & 0.1s & 0.1s & 0.1s & 0.1s & 0.1s \\
 $n = 20$   & 3s & 1s & 0.3s  & 0.2s & 0.1s \\
 $n = 30$   & 17s & 5s & 2s & 0.7s & 0.2s \\
 $n = 40$  & 55s & 16s & 4s & 2s & 0.5s \\
 $n = 50$  & 112s & 33s & 8s & 4s & 1s \\
\midrule
\emph{Instances} & \multicolumn{5}{l}{\emph{Average solution time of network flow model, in seconds}} \\ 
 $n = 10$  & 0.1s & 0.1s & 0.1s & 0.1s & 0.1s \\
 $n = 20$   & 1s & 0.5s  & 0.3s & 0.2s & 0.1s \\
 $n = 30$  & 87s & 11s & 3s & 1s & 0.3s \\
 $n = 40$  & 534s & 45s & 8s & 3s & 1s \\
 $n = 50$  & $>$3600s & 117s & 22s & 8s & 2s \\
\bottomrule
\end{tabular}}
\caption{The average build time and solution time of relaxed decision diagrams under different $Q$ values. Build time includes time to generate and reduce a diagram.}\label{table:RDD_times}
\end{table}

Next, we examine the quality of the solution generated by these approximate models.

\subsubsection{Quality of model solutions.} For each solution, we calculate the model-based optimality gap as well as the true optimality gap where possible (see Section \ref{subsec:evaluation_oracle} for details). Table \ref{table:RDD_optgap} summarizes the values of these optimality gaps. Note that only the model-based optimality gap is shown for the $n = 50$ instances, since the exact models (i.e., where $Q = 0$) could not be solved within the time limit. We also note that since the formulation and solution times of the exact models are generally negligible for $n = 10$ and $n = 20$ (see Table \ref{table:RDD_times}), we focus our discussion on the instances with $n = 30, 40$ and $50$, where tractable approximations become more critical.

\begin{table}[h]
\centering
\setlength{\tabcolsep}{12pt}
{\renewcommand{\arraystretch}{1.2}%
\begin{tabular}{c|rrrrr}
 \toprule
 Instances & $Q = 0$ & $Q = 1$ & $Q = 3$ & $Q = 5$ & $Q = 10$ \\
\midrule 
 $n = 10$ & 0\%  & 1.2\% (1.5\%) & 3.0\% (4.3\%) & 2.6\% (4.6\%) & 4.4\% (7.6\%)  \\ 
 $n = 20$ & 0\% & 1.2\% (1.5\%) & 1.5\% (2.1\%) & 1.0\% (1.7\%) & 1.1\% (2.2\%)  \\ 
 $n = 30$ & 0\% & 0.3\% (0.5\%) & 0.4\% (0.7\%) & 0.5\% (0.8\%) & 0.9\% (1.4\%)  \\ 
 $n = 40$ & 0\% & 0.1\% (0.2\%) & 0.2\% (0.3\%) & 0.2\% (0.4\%) & 0.2\% (0.5\%) \\ 
 $n = 50$ & - & - (0.08\%) & - (0.2\%) & - (0.3\%) & - (0.5\%)  \\
\bottomrule
\end{tabular}}
\vspace{0.1cm}
\caption{The average true optimality gap and model-based optimality gap of the solution, the latter of which is shown in parentheses.}\label{table:RDD_optgap}
\end{table}

The main takeaway from Table \ref{table:RDD_optgap} is that both the true and the model-based optimality gaps remain very small despite the large reduction in solution times (as shown in Table \ref{table:RDD_times}). For example, when $n = 40$ and $Q = 5$, the average optimality gap of solutions is less than 0.5\%, despite the models taking two orders of magnitude less time to formulate and solve compared to the exact model (i.e., an average of 5 seconds versus 589 seconds; see Table \ref{table:RDD_times}). Similarly, when $n = 50$, the average model-based optimality gap is $0.08\%$ when $Q = 1$ and $0.5\%$ when $Q = 10$, the latter of which requires an average solution time that is three orders of magnitude less than the exact model (i.e., 3 seconds versus $3600+$ seconds; see Table \ref{table:RDD_times}). These observations highlight that the approximate network flow models are able to independently generate both (i) near-optimal solutions \emph{and} (ii) almost-tight dual bounds, as both of these conditions must be met to observe small model-based optimality gaps. Finally, as a side note, we point out that for any fixed value of $Q$, the average optimality gap decreases as $n$ increases. This is because in our problem setting, by holding $Q$ constant, the relative degree of the approximation decreases as the problem size increases.

\begin{figure}[h]
    \captionsetup[subfigure]{skip=3pt}
    \centering
    \begin{subfigure}{0.49\textwidth}
    \centering
    \includegraphics[width = 7.9cm]{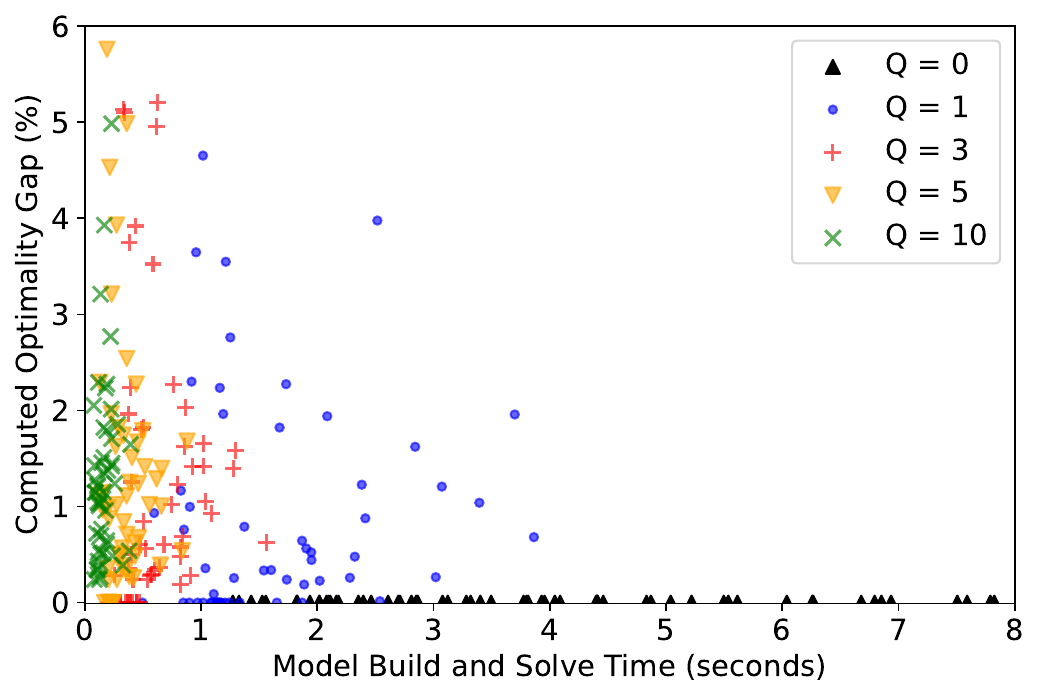}
    \captionsetup{font={small}}
    \caption{n = 20}
    \end{subfigure}
    \begin{subfigure}{0.48\textwidth}
    \centering
    \includegraphics[width = 7.9cm]{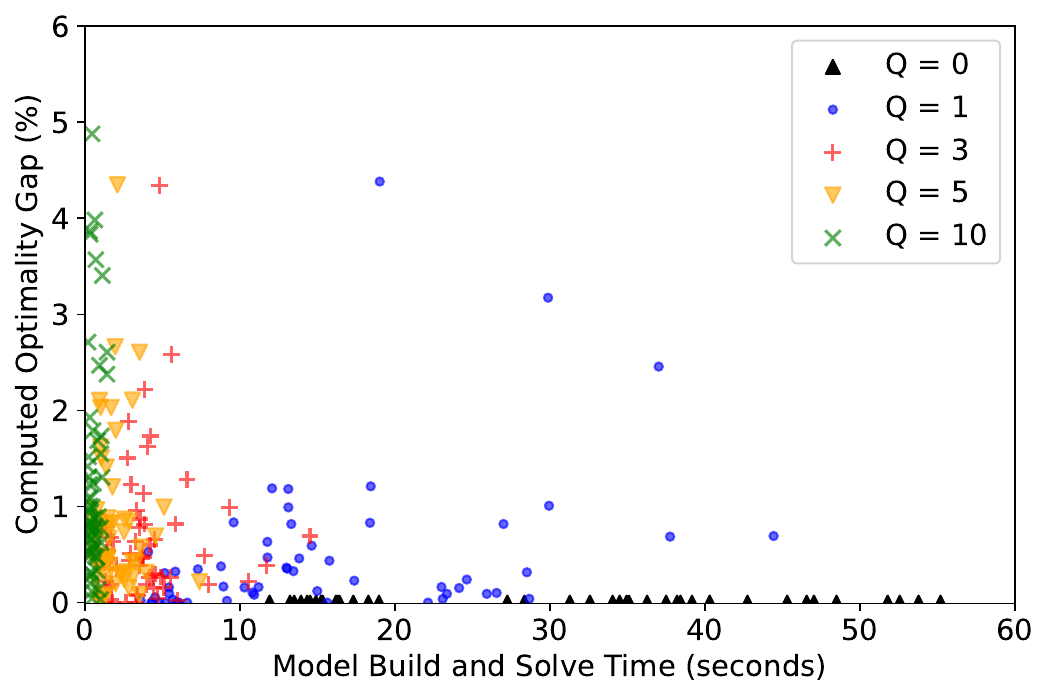}
    \captionsetup{font={small}}
    \caption{n = 30}
    \end{subfigure}
    
    \begin{subfigure}{0.49\textwidth}
    \centering
    \includegraphics[width = 7.9cm]{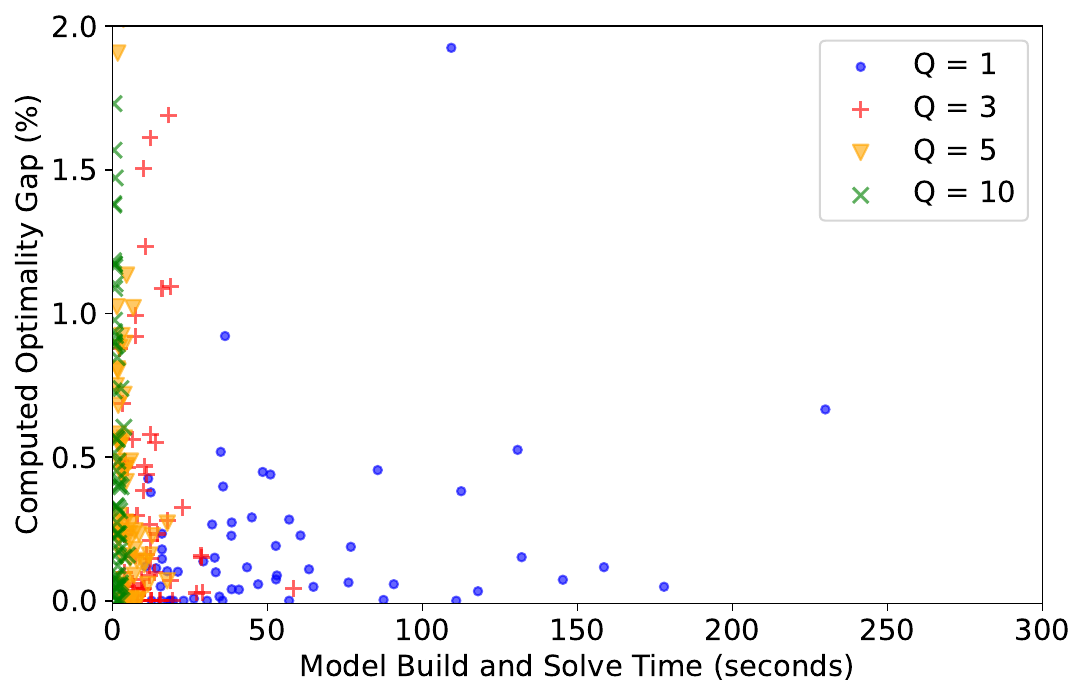}
    \captionsetup{font={small}}
    \caption{n = 40}
    \end{subfigure}
    \begin{subfigure}{0.48\textwidth}
    \centering
    \includegraphics[width = 7.9cm]{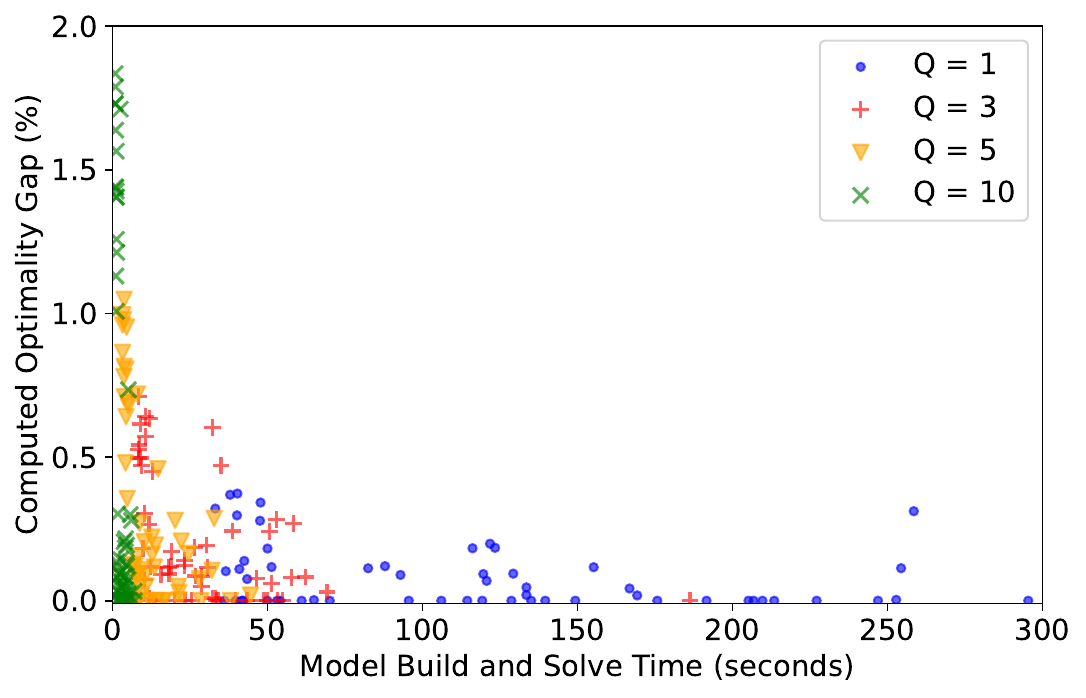}
    \captionsetup{font={small}}
    \caption{n = 50}
    \end{subfigure}
    \vspace{0.2cm}
    \caption{The model-based optimality gap as a function of the total time to formulate and solve the network flow models. Note that the data points roughly follow an exponential decay curve, that is, an incremental increase in solution time can result in a large decrease in optimality gap, when solution times are small. }\label{fig:CB_optgap}
\end{figure}

Figure \ref{fig:CB_optgap} provides a more nuanced illustration of the model-based optimality gap for each instance as a function of solution time. The main observation is that across all problem sizes, we observe an exponential decay in the optimality gap as a function of solution time. Specifically, the results show that when solution times are small, a slight increase in model complexity and solution time (i.e., by decreasing the value of $Q$ slightly) can result in a large decrease in optimality gap. Put differently, a slight approximation of the exact model can drastically reduce solution times while sacrificing very little in terms of the quality of solutions and dual bounds. This is most evident in the $n=40$ ($n=50$) instances, where reducing the average solution time from 589 seconds (3600+ seconds) to 1.5 seconds (3 seconds) results in an average suboptimality loss of 0.2\% (at most 0.5\%). These insights highlight that our approximate network flow models are highly effective and scalable.

\subsubsection{Brief comparison with alternative methods.}

We briefly comment on the performance of the $K$-adaptability model \citep{hanasusanto2015k}, which is a popular approximation method that has also served as a benchmark for the same capital budgeting instances (i.e., see \cite{arslan2022decomposition} and \cite{dumouchelle2023neur2ro}). The $K$-adaptability model generates a feasible solution and a primal bound; for reference, we provide the complete formulation in Section \ref{ECsub:KAdaptCapital}. Since there is strong evidence from previous literature that the model can be challenging to solve even for small problem sizes and small $K$, we, for practical reasons, select a sub-sample of 100 instances and impose a maximum time limit of 1200 seconds for each instance. Our sub-sample consists of 20 instances for each value of $n = \{10, 20, 30, 40, 50\}$, which we select simply by choosing every third instance in the GitHub repository from \cite{arslan2022decomposition}.

Table \ref{table:K-adaptCB} highlights the solution times of the $K$-adaptability model for $K = 2,3,4$. First, note that the solution times for the $K$-adaptability model increase exponentially as we increase either the value of $K$ or the size of the instances.  For example, even when $K = 2$, many instances where $n \geq 30$ cannot be solved within the 1200-second time limit. As another example, when $n=20$, the average solution time goes from 13 seconds to over 675 seconds when $K$ is increased from 2 to 4. When $n = 40$ and $n = 50$, 15 out of 20 instances could not be solved within 1200 seconds for any $K$, while 5 of 20 took less than one second (resulting in the consistent average of 900 seconds). Finally, many instances of $K$-adaptability had large optimality gaps (e.g., $> 10\%$) when the 1200-second time limit was reached.

\begin{table}[h]
\centering
\setlength{\tabcolsep}{12pt}
{\renewcommand{\arraystretch}{1.2}%
\begin{tabular}{c|rrr|rrr}
 \toprule
Instances & $Q = 0$ & $Q = 5$ & $Q = 10$ & \ \ \  $ K = 2$ & $K = 3$ & $ K = 4$ \\
\midrule 
n = 10  & $0.1$s & $0.1$s & $0.1$s & $0.1$s & $0.2$s & $1.2$s \\ 
n = 20  & $4$s & $0.3$s & $0.1$s & $13$s & $>171$s & $>675$s  \\ 
n = 30  & $176$s & $2$s & $0.5$s & $>627$s & $>755$s & $>900$s  \\ 
n = 40  & - & $4$s & $2$s & $>900$s & $>900$s & $>900$s  \\ 
n = 50  & - & $10$s & $3$s & $>900$s & $>900$s & $>900$s  \\
\bottomrule
\end{tabular}}
\vspace{0.1cm}
\caption{A comparison of average solution times of approximation models. For each entry, the symbol $>$ is used to denote any average that is taken when there exists at least 1 instance that exceeds the 1200-second threshold. The total time for the network flow models include the model formulation time (i.e., diagram generation time).}\label{table:K-adaptCB}
\end{table}

In comparison to the $K$-adaptability models, our network flow models can be solved much more efficiently. For example, the average solution time for instances where $n = 50$ and $Q = 5$ is $10$ seconds, and, as we highlighted in the previous subsection, the generated first-stage solutions are within $0.2\%$ of optimality. As another example, it takes an average of 4 seconds to solve the $n = 20$ instances exactly, while it takes more than 675 seconds to solve the $K$-adaptability model with $K=4$ (which still does not generate optimal solutions for all instances). Finally, and perhaps most importantly, our network flow models also simultaneously generate high-quality dual bounds, which allow us to evaluate the quality of \emph{any} feasible solution, including those that are generated independently by heuristics or other approximation models like $K$-adaptability.

In summary, our models provide a flexible framework for generating both high-quality solutions and dual bounds in significantly less time. Compared to $K$-adaptability, the complexity of our model can also be tuned much more precisely. For example, increasing the value of $Q$ gradually increases the solution time of the model, whereas solution times increase exponentially with small changes in the value of $K$. Related discussion was also presented in Sections \ref{subsubsec:distance-based} and \ref{subsec:primal_bound}.

Lastly, we remark that while our discussion focuses on the $K$-adaptability model for reasons previously mentioned (e.g., popularity and ease of implementation in standard solvers), we also compare our results to the solutions times of the exact branch-and-price algorithm presented in \cite{arslan2022decomposition}. We find that the solution times of many exact and/or near-exact reformulations do not exceed those reported for the branch-and-price algorithm, but more importantly, our approach can generate approximate reformulations that are significantly faster to solve while sacrificing little in terms of solution quality. We refer the reader to Section \ref{ECsub:BPresults} of the Electronic Companion for details.

\subsection{Robust assignment problems}\label{subsec:assignment}

In the previous subsection, we examined exact and approximate single-network flow models in the context of the capital budgeting problem, which has a single linking constraint between first-stage decisions and each second-stage decision. In this subsection, we examine a problem setting in which numerous such constraints exist. Specifically, we examine robust assignment problems and focus on the use of multi-network flow models, which are discussed in Section \ref{subsec:dual_bounds}.

\begin{figure}
\centering
\begin{tikzpicture}[>=stealth]
\tikzset{square/.style={rectangle, draw, minimum width=0.6cm, minimum height=0.6cm}}

\def\nodeDist{1.5cm}
\def\rightlabelDist{-1.3cm}
\def\leftlabelDist{0.45cm}

\node[circle, draw, minimum size=0.6cm] (leftL) at (0,1*\nodeDist) {\vspace{2pt}$L$\vspace{2pt}};
\node (leftDots) at (0,2*\nodeDist) {\vspace{2pt}$\vdots$\vspace{2pt}};
\node[circle, draw, minimum size=0.6cm] (left2) at (0,3*\nodeDist) {\vspace{2pt}$2$\vspace{2pt}};
\node[circle, draw, minimum size=0.6cm] (left1) at (0,4*\nodeDist) {\vspace{2pt}$1$\vspace{2pt}};

\node[square, draw] (rightM) at (4,1.5*\nodeDist) {\vspace{2pt}$M$\vspace{2pt}};
\node[] (rightDots) at (4,2.5*\nodeDist) {\vspace{2pt}$\vdots$\vspace{2pt}};
\node[square, draw] (right1) at (4,3.5*\nodeDist) {\vspace{2pt}$1$\vspace{2pt}};

\draw[->] (left1) -- (right1) node[midway, above] {\small$\xi_{1,1}$};
\draw[->] (left1) -- (rightDots) {};
\draw[->] (left1) -- (rightM) {};
\draw[->] (left2) -- (right1) {};
\draw[->] (leftL) -- (rightDots) {};
\draw[->] (leftL) -- (rightM) node[midway, above] {\small$\xi_{L,M}$};

\node[left] at (left1) {\hspace{\rightlabelDist}\small $a_1$};
\node[left] at (left2) {\hspace{\rightlabelDist}\small $a_2$};
\node[left] at (leftL) {\hspace{\rightlabelDist}\small $a_L$};

\node[right] at (right1) {\hspace{\leftlabelDist}\small $b_1$};
\node[right] at (rightM) {\hspace{\leftlabelDist}\small $b_M$};

\end{tikzpicture}
\caption{A robust assignment problem. Agents $1, \ldots, L$ are associated with weight $a_1, \ldots, a_L$, tasks $1, \ldots, M$ are associated with capacity $b_1, \ldots, b_M$, and agent-task pairings $(\ell,m)$ generate reward $\xi_{\ell m}$.}
\label{fig:assignment}
\end{figure}
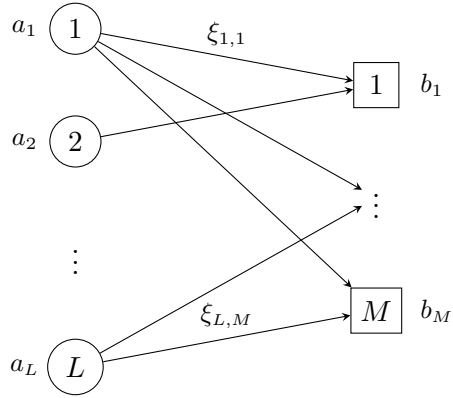

Assignment problems encompass numerous decision-making tasks that span many applications. A standard assignment problem can be modeled as a bipartite graph $(\mL, \mM, \mS)$ of agents $\mL = \{1, \ldots, L\}$, tasks $\mM = \{1, \ldots, M\}$ and directed links $\mS$. Let $\mL(m) \subseteq \mL$ denote the subset of agents for which there exists a directed link to task $m \in \mM$, and similarly, let $\mM(\ell) \subseteq \mM$ denote the subset of tasks which can be assigned to agent $\ell \in \mL$. Depending on the context, agents can represent people, products, resources, funding, or jobs, whereas tasks can represent locations, facilities, projects or machines. Each agent $\ell \in \mL$ is associated with weight $a_\ell \geq 0$ 
and each task $m$ has capacity $b_m \geq 0$, 
while the reward $\xi_{\ell,m} \geq 0$ (e.g., match quality) of a specific agent-task pairing is unknown (e.g., varies day-to-day or must be estimated from data). Figure \ref{fig:assignment} illustrates this assignment problem.

Given this setup, a standard robust assignment problem can be formulated as 
\begin{alignat*}{2}
\underset{\bx \in \mX}{\text{max}} \ \ \underset{\bxi \in \Xi}{\text{min}} \ \ \bxi^\top \bx \qquad \text{with} \qquad \mX = \Big \{\, \bx \in \{0,1\}^{|\mS|} ~:~ & \sum_{\ell \in \mL(m)} a_\ell x_{\ell,m} \leq b_m, \quad && \forall m \in \mM, \\
& \sum_{m \in \mM(\ell)} x_{\ell,m} \leq 1, \quad && \forall \ell \in \mL \ \Big \}.  
\end{alignat*}

Consider a setting where a planner seeks to introduce a limited degree of flexibility in the assignment planning and decision-making process. Specifically, suppose that the planner wants to pre-identify a set of potential agent-task pairings from which assignments can be made once rewards become known. Practically, this may arise when there is a desire to train people for specific tasks, inform individuals about potential assignments, or inspect match quality before making a final decision. We assume that we have a cardinality constraint limiting the number of pre-identified pairings, and formulate the adaptive robust assignment problem as
\begin{subequations}\label{model:2S-assignment}
\begin{alignat}{2}
\underset{\bx \in \mX}{\text{max}} \ \ \underset{\bxi \in \Xi}{\text{min}} \ \ \underset{\by}{\text{max}} \quad & \ \bxi^\top \by\\
\text{s.t.} \quad & \sum_{\ell \in \mL(m)} a_\ell y_{\ell,m} \leq b_m, \quad && \forall m \in \mM \label{cons:2S_assign1} \\[0.1cm]
& \sum_{m \in \mM(\ell)} y_{\ell,m} \leq 1, \quad && \forall \ell \in \mL \label{cons:2S_assign2} \\[0.1cm]
& \ 
\by \in \{0,1\}^{|\mS|} \label{cons:2S_binary} \\
& \ \by \leq \bx 
\end{alignat}
\end{subequations}
where $\mX = \left\{ \bx \in \{0,1\}^{|\mS|} ~:~ \norm{\bx}_1 \leq \beta \cdot |\mS| \right\}$. 
In our experiments, we consider $\beta$ values of 0.5, 0.6, 0.7, 0.8, and 0.9, which correspond to the ability to pre-identify $50\%$, $60\%$, $70\%$, $80\%$, and $90\%$ of pairings in the first stage, which can then be used to form a final assignment in the second stage. 

\subsubsection{Experimental setup.} We consider assignment problems of five sizes where $(L,M)$ is equal to (20,2), (20,3), (20,4), (25,5) and (25,8). For each problem size, we randomly generate 10 instances. Each instance is characterized by (i) a randomly generated a bipartite graph with 50\% sparsity (such that $|\mS| = 20, 30, 40, 63$ and $100$, respectively) and (ii) randomly generated vectors of agent weights $\ba$ and task capacities $\bb$, where each element of $\ba$ is a random integer from $[1, 10]$ and each element of $\bb$ is an random integer from $[2\cdot\text{max}(\ba), \frac{1}{m}||\ba||_1]$.

For each instance, we generate the uncertainty set as follows. Let $\bxi^0$ denote a vector representing the nominal reward of each link. With slight abuse of notation, each element $\xi^0_i$ is a random number generated from $[0.5a, a]$ where $a$ is the weight of the agent associated with this link. 
Then, the uncertainty set is defined on the percentage of deviation from $\bxi^0$, namely,
\begin{align*}
\Xi = \left\{\bxi \; \bigg | \; \sum_{i = 1}^{|\mS|} \big |\xi_i/\xi^0_i - 1 \big| \leq 0.1|\mS|, \ \big |\xi_i/\xi^0_i - 1 \big | \leq 0.5 \ \ \forall i \in \{1, \ldots, |\mS|\} \right\}.  
\end{align*}

We consider both an exact network flow formulation and a multi-network flow approximation for solving the adaptive assignment problems. These two formulations are denoted as Exact NF and Multi NF in the tables and figures, and we generate them as follows:
\begin{itemize}
\item\textbf{Exact NF.} Model \eqref{model:2S-assignment} can be represented as an exact network flow model by reformulating constraints \eqref{cons:2S_assign1} -- \eqref{cons:2S_binary} based on a straightforward extension of the procedure outlined in Example \ref{ex:BDD_construction}. Specifically, the state $\bS$ in the recursive formulation of \eqref{cons:2S_assign1} -- \eqref{cons:2S_binary} is a vector of size $L+M$ (rather than scalars) that captures the amount of capacity remaining in \eqref{cons:2S_assign1} -- \eqref{cons:2S_assign2} based on previous decisions of $y_1, \ldots, y_i$. Generating this recursive formulation (i.e., the decision diagram) is straightforward, as shown in Section \ref{appsub:RAP_recursion} of the Electronic Companion.

\item \textbf{Multi NF.} To form the multi-network flow approximation, we generate an exact decision diagram for each knapsack constraint in \eqref{cons:2S_assign1} along with binary domains in \eqref{cons:2S_binary} and integrate the corresponding network flow constraints. We leave constraints \eqref{cons:2S_assign2} intact since they satisfy the integral polyhedron property (as an implication of Proposition \ref{propn:constraint_decoupling} in Section \ref{subsec:dual_bounds}). This model generates a feasible first-stage solution and a dual bound. 

\end{itemize}

The Multi NF model defined above can be considered as one of the most natural or straightforward approximations of the ARBO problem, and we will show that it performs well in the numerical experiments. Nonetheless, we remark that there are many possible variations of this model which we could use to tighten or loosen the approximation. For example, instead of generating one exact decision diagram for each knapsack constraint, we could generate one for each pair of constraints, which would tighten the approximation quality but potentially lead to larger formulations. On the other hand, we could also generate an approximate decision diagram for each knapsack constraint (like in Section \ref{subsec:capital_budgeting}), which leads to a worse approximation but could significantly reduce the size and solution time of the model. In summary, various Multi NF models could be defined that will trade off between lower solution times and better solution quality.

Finally, for each network flow model, we enforce a one-hour time limit to build and reduce each diagram. Similar to Section \ref{subsec:capital_budgeting}, we also consider the $K$-adaptability model with $K=3$ and $K=4$ as a point of reference; the complete MILP formulation can be found in Section \ref{appsub:RAP_Kadapt} of the Electronic Companion. For all models, we enforce a 30-minute time limit for the solver.

\begin{table}[h]
\centering
\setlength{\tabcolsep}{6pt}
{\renewcommand{\arraystretch}{1.1}%
\begin{tabular}{ll rrrrrr}
 \toprule
Model & BDD attributes & (20,2) & (20,3) & (20,4) & (25,5) & (25, 8) \\[0.1cm]
 \midrule
  \multirow{4}{*}{Exact NF} & avg. build time & 0.5s & 67s & \ \ - & \ \ - & \ \ - \\
& avg. reduction time & 1.4s & 243s & \ \ - & \ \ - & \ \ -  \\
& avg. \# of arcs (unreduced) & 11994 & 199856 & \ \ - & \ \ - & \ \ - \\
& avg. \# of arcs (reduced) & 498 & 5506 & \ \ - & \ \ - & \ \ - \\
 \midrule
\multirow{4}{*}{Multi NF} & avg. build time & 0.1s & 0.1s & 0.1s & 0.1s & 0.1s \\
& avg. reduction time & 0.1s & 0.1s & 0.1s & 0.2s & 0.4s \\
& avg. \# of arcs (unreduced) & 1262 & 2577 & 3911 & 8855 & 20365 \\
& avg. \# of arcs (reduced) & 287 & 884 & 1385 & 3669 & 7742\\
\bottomrule
\end{tabular}}
\caption{The build time and attributes of the decision diagrams. Values for Multi NF are given as a summation over all the diagrams used to generate the multi-network formulation for a particular instance. 
}\label{table:BDD_sizes}
\end{table}

\subsubsection{Computational results.}

We first give an overview of the size and formulation time of the network flow models. Then, we examine the solution time, solution quality, and dual bounds generated by the models.

The attributes of the decision diagrams underlying the network flow models are shown in Table \ref{table:BDD_sizes}. Specifically, Table \ref{table:BDD_sizes} highlights the time it takes to generate and reduce the diagrams as well as the sizes of the diagrams. We note that the size and formulation time of the exact decision diagrams increases rapidly with the size of the assignment problem. Specifically, it takes an average of 2 seconds for the (20,2) instances, 310 seconds for the (20,3) instances, and more than one hour for all the other instances. Since the set of all feasible recourse decisions must be represented within a single decision diagram, the size of this diagram may grow exponentially with the number of constraints used to define the feasible space. In fact, in the robust assignment problem, the main computational bottleneck problem is in generating the exact decision diagram rather than solving the corresponding network flow formulation. This will become even more clear in the next paragraph. In contrast, the average size and formulation time of the multi-network models is far smaller and generally negligible ($<1$ second). Since this model is simply a collection of individual decision diagrams that are each exact reformulations of only one constraint in \eqref{cons:2S_assign1}, the size and formulation time scales linearly in the number of constraints in \eqref{cons:2S_assign1}. Note that for Multi NF models, the generation of individual diagrams for each constraint could be parallelized.

Table \ref{table:exact_times} highlights the average solution time of the various models over different problem instances. Note that despite the exact decision diagram taking significant time to generate for the (20,3) instances, it only takes 1 second to solve the exact network flow formulation. As for the Multi NF models, they can be solved relatively efficiently even as the size of the assignment problem grows larger. We briefly point out that the solution times for the $K$-adaptability model with $K = 3$ and $K=4$ are orders of magnitude larger than the network flow models.

\begin{table}[h]
\centering
\setlength{\tabcolsep}{5.8pt}
\begin{tabularx}{0.92\textwidth} { l  
   rrrrrrrrrr}
 \toprule
 Model & \multicolumn{2}{c}{(20,2)} & \multicolumn{2}{c}{(20,3)} & \multicolumn{2}{c}{(20,4)} & \multicolumn{2}{c}{(25,5)} & \multicolumn{2}{c}{(25,8)} \\
 \cmidrule(lr){1-1}
 \cmidrule(lr){2-3}
 \cmidrule(lr){4-5}
 \cmidrule(lr){6-7}
 \cmidrule(lr){8-9}
 \cmidrule(lr){10-11}
 Exact NF & 0.1s & (50) & 1.1s &  (50) & - &  & - &  & - &   \\
Multi NF & 0.1s & (50) & 0.3s & (50) & 1.8s & (50) & 81s & (50) & $>$633s & (33)  \\[0.15cm]
\midrule\\[-0.25cm]
3-Adapt & 1.2s & (50) & 48s & (50) & $>324$s & (46) & $>1634$s & (5) & $>1800$s & (0) \\
4-Adapt & 26s & (50) & $>342$s & (46) & $>1134$s & (24) & $>1749$s & (3) & $>1800$s & (0)  \\
 \bottomrule
\end{tabularx}
\caption{Average solution times across the different models. The value in the parentheses denotes the number of instances, out of 50, that could be solved within the 30-minute time limit given to each instance.}\label{table:exact_times}
\end{table}

\begin{figure}[h]
    \centering
    \includegraphics[width = 10cm]{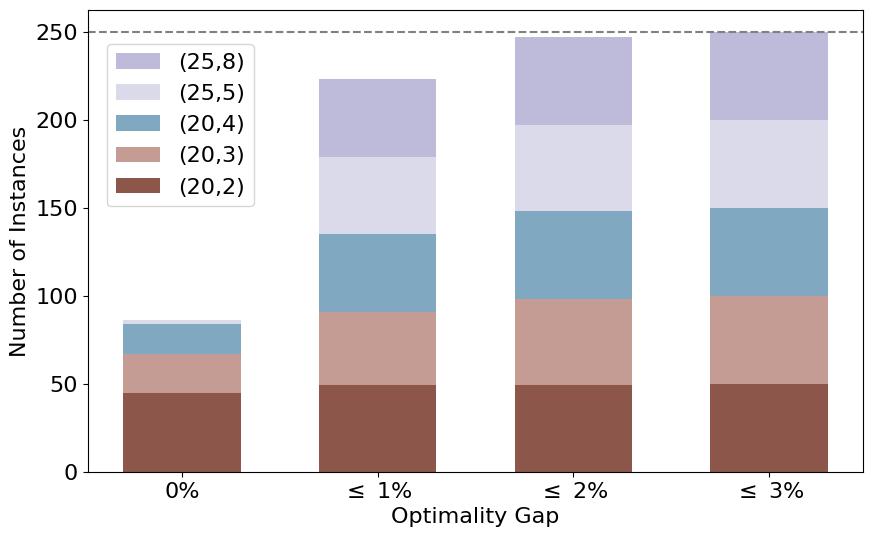}
    \caption{Model-based optimality gap of solutions of the Multi NF model.}\label{fig:RAP_optgaps}
\end{figure}

Despite having significantly lower solution times, the Multi NF models are still able to generate solutions that have very low model-based optimality gaps, as shown in Figure \ref{fig:RAP_optgaps}. Note that, as discussed in Section \ref{subsec:capital_budgeting}, these values are upper bounds on the true optimality gap of the solution, meaning that the true quality of the solution could be even better than expected. First, we note that 86 out of the 250 problem instances had a model-based optimality gap of 0\%. In these instances, the Multi NF model found the optimal solution \emph{and} verified its optimality. Most of these cases pertained to problem instances that were smaller in size. Second, 223 out of 250 instances had a model-based optimality gap that was less than 1\%, 247 out of 250 had a gap less than 2\%, and no solution had an optimality gap that was greater than 3\%. These statistics include solutions of the Multi NF instances that could not be solved within the time limit (see Table \ref{table:exact_times}). 

The main takeaway of these results is that the multi-network flow model can efficiently generate near-optimal first-stage solutions \emph{and} high-quality dual bounds across all problem instances considered. This is consistent with the observations from Section \ref{subsec:capital_budgeting}, which highlighted the same findings using approximate single-network flow models.

\subsection{Summary of numerical results}

We briefly summarize the main takeaways of the numerical section. First, we show that exact network flow formulations can be tractably formulated and solved for smaller instances of ARBO problems. Second, for larger instances, the approximate network flow formulations are tractable and simultaneously generate (i) high-quality solutions and (ii) high-quality dual bounds. Finally, the size and solution times of the approximate network flow models can often be reduced drastically while sacrificing very little in terms of the quality of solutions and dual bounds.

\section{Conclusion}\label{sec:conclusion}

In this paper, we examine adaptive robust binary optimization problems with objective uncertainty. We leverage ideas from the decision diagram community to reformulate and approximate our adaptive problems as single-stage constrained network flow models. We outline methods to generate network models where the size and quality of the models can be easily controlled through a user-specified parameter. Our models are also easy to implement and solve using standard MILP solvers. Through an extensive set of computational experiments, we show that these models can efficiently generate both high-quality solutions and high-quality dual bounds.

By forming this connection between adaptive robust optimization and decision diagrams, our framework can also take advantage of independent research contributions that emerge from the latter research community. Specifically, these developments may create opportunities to extend our ideas to more general problem settings. For example, recent literature has considered the use of decision diagrams to solve deterministic integer and/or nonlinear optimization problems \citep{castro2022decision}, and similar ideas could potentially be integrated into our framework to reformulate adaptive robust problems with general integer recourse.

\bibliographystyle{plainnat}
\bibliography{SROptimization}

\begin{thebibliography}{36}
\providecommand{\natexlab}[1]{#1}
\providecommand{\url}[1]{\texttt{#1}}
\expandafter\ifx\csname urlstyle\endcsname\relax
  \providecommand{\doi}[1]{doi: #1}\else
  \providecommand{\doi}{doi: \begingroup \urlstyle{rm}\Url}\fi

\bibitem[{\'A}lvarez-Miranda et~al.(2015){\'A}lvarez-Miranda, Fern{\'a}ndez, and Ljubi{\'c}]{alvarez2015recoverable}
Eduardo {\'A}lvarez-Miranda, Elena Fern{\'a}ndez, and Ivana Ljubi{\'c}.
\newblock The recoverable robust facility location problem.
\newblock \emph{Transportation Research Part B: Methodological}, 79:\penalty0 93--120, 2015.

\bibitem[Arslan and Detienne(2022)]{arslan2022decomposition}
Ay{\c{s}}e~N Arslan and Boris Detienne.
\newblock Decomposition-based approaches for a class of two-stage robust binary optimization problems.
\newblock \emph{INFORMS Journal on Computing}, 34\penalty0 (2):\penalty0 857--871, 2022.

\bibitem[Bayram et~al.(2022)Bayram, Baloch, Gzara, and Elhedhli]{bayram2022optimal}
Vedat Bayram, Gohram Baloch, Fatma Gzara, and Samir Elhedhli.
\newblock Optimal order batching in warehouse management: A data-driven robust approach.
\newblock \emph{INFORMS Journal on Optimization}, 2022.

\bibitem[Bergman et~al.(2016)Bergman, Cire, Van~Hoeve, and Hooker]{bergman2016decision}
David Bergman, Andre~A Cire, Willem-Jan Van~Hoeve, and John Hooker.
\newblock \emph{Decision diagrams for optimization}, volume~1.
\newblock Springer, 2016.

\bibitem[Bergman et~al.(2022)Bergman, Bodur, Cardonha, and Cire]{bergman2022network}
David Bergman, Merve Bodur, Carlos Cardonha, and Andre~A Cire.
\newblock Network models for multiobjective discrete optimization.
\newblock \emph{INFORMS Journal on Computing}, 34\penalty0 (2):\penalty0 990--1005, 2022.

\bibitem[Bertsimas and Dunning(2016)]{bertsimas2016multistage}
Dimitris Bertsimas and Iain Dunning.
\newblock Multistage robust mixed-integer optimization with adaptive partitions.
\newblock \emph{Operations Research}, 64\penalty0 (4):\penalty0 980--998, 2016.

\bibitem[Bertsimas and Georghiou(2015)]{bertsimas2015design}
Dimitris Bertsimas and Angelos Georghiou.
\newblock Design of near optimal decision rules in multistage adaptive mixed-integer optimization.
\newblock \emph{Operations Research}, 63\penalty0 (3):\penalty0 610--627, 2015.

\bibitem[Bertsimas et~al.(2013)Bertsimas, Nasrabadi, and Stiller]{bertsimas2013robust}
Dimitris Bertsimas, Ebrahim Nasrabadi, and Sebastian Stiller.
\newblock Robust and adaptive network flows.
\newblock \emph{Operations Research}, 61\penalty0 (5):\penalty0 1218--1242, 2013.

\bibitem[Bryant(1992)]{bryant1992symbolic}
Randal~E Bryant.
\newblock Symbolic boolean manipulation with ordered binary-decision diagrams.
\newblock \emph{ACM Computing Surveys (CSUR)}, 24\penalty0 (3):\penalty0 293--318, 1992.

\bibitem[Buchheim and Kurtz(2018)]{buchheim2018robust}
Christoph Buchheim and Jannis Kurtz.
\newblock Robust combinatorial optimization under convex and discrete cost uncertainty.
\newblock \emph{EURO Journal on Computational Optimization}, 6\penalty0 (3):\penalty0 211--238, 2018.

\bibitem[Castro et~al.(2022)Castro, Cire, and Beck]{castro2022decision}
Margarita~P Castro, Andre~A Cire, and J~Christopher Beck.
\newblock Decision diagrams for discrete optimization: {A} survey of recent advances.
\newblock \emph{INFORMS Journal on Computing}, 34\penalty0 (4):\penalty0 2271--2295, 2022.

\bibitem[Cire et~al.(2019)Cire, Diamant, Yunes, and Carrasco]{cire2019network}
Andre~A Cire, Adam Diamant, Tallys Yunes, and Alejandro Carrasco.
\newblock A network-based formulation for scheduling clinical rotations.
\newblock \emph{Production and Operations Management}, 28\penalty0 (5):\penalty0 1186--1205, 2019.

\bibitem[Da{\c{s}} et~al.(2020)Da{\c{s}}, Gzara, and St{\"u}tzle]{dacs2020review}
G{\"u}lesin~Sena Da{\c{s}}, Fatma Gzara, and Thomas St{\"u}tzle.
\newblock A review on airport gate assignment problems: Single versus multi objective approaches.
\newblock \emph{Omega}, 92:\penalty0 102146, 2020.

\bibitem[de~Lima et~al.(2022)de~Lima, Alves, Clautiaux, Iori, and de~Carvalho]{de2022arc}
Vin{\'\i}cius~L de~Lima, Cl{\'a}udio Alves, Fran{\c{c}}ois Clautiaux, Manuel Iori, and Jos{\'e} M~Val{\'e}rio de~Carvalho.
\newblock Arc flow formulations based on dynamic programming: Theoretical foundations and applications.
\newblock \emph{European Journal of Operational Research}, 296\penalty0 (1):\penalty0 3--21, 2022.

\bibitem[Dumouchelle et~al.(2023)Dumouchelle, Julien, Kurtz, and Khalil]{dumouchelle2023neur2ro}
Justin Dumouchelle, Esther Julien, Jannis Kurtz, and Elias~B Khalil.
\newblock Neur2ro: Neural two-stage robust optimization.
\newblock \emph{arXiv preprint arXiv:2310.04345}, 2023.

\bibitem[Eufinger et~al.(2020)Eufinger, Kurtz, Buchheim, and Clausen]{eufinger2020robust}
Lars Eufinger, Jannis Kurtz, Christoph Buchheim, and Uwe Clausen.
\newblock A robust approach to the capacitated vehicle routing problem with uncertain costs.
\newblock \emph{INFORMS Journal on Optimization}, 2\penalty0 (2):\penalty0 79--95, 2020.

\bibitem[Gorissen et~al.(2015)Gorissen, Yan{\i}ko{\u{g}}lu, and den Hertog]{gorissen2015practical}
Bram~L Gorissen, {\.I}hsan Yan{\i}ko{\u{g}}lu, and Dick den Hertog.
\newblock A practical guide to robust optimization.
\newblock \emph{Omega}, 53:\penalty0 124--137, 2015.

\bibitem[Guo et~al.(2021)Guo, Bodur, Aleman, and Urbach]{guo2021logic}
Cheng Guo, Merve Bodur, Dionne~M Aleman, and David~R Urbach.
\newblock Logic-based {B}enders decomposition and binary decision diagram based approaches for stochastic distributed operating room scheduling.
\newblock \emph{INFORMS Journal on Computing}, 33\penalty0 (4):\penalty0 1551--1569, 2021.

\bibitem[Hanasusanto et~al.(2015)Hanasusanto, Kuhn, and Wiesemann]{hanasusanto2015k}
Grani~A Hanasusanto, Daniel Kuhn, and Wolfram Wiesemann.
\newblock K-adaptability in two-stage robust binary programming.
\newblock \emph{Operations Research}, 63\penalty0 (4):\penalty0 877--891, 2015.

\bibitem[Hooker(2013)]{hooker2013decision}
John~N Hooker.
\newblock Decision diagrams and dynamic programming.
\newblock In \emph{International Conference on Integration of Constraint Programming, Artificial Intelligence, and Operations Research}, pages 94--110. Springer, 2013.

\bibitem[K{\"a}mmerling and Kurtz(2020)]{kammerling2020oracle}
Nicolas K{\"a}mmerling and Jannis Kurtz.
\newblock Oracle-based algorithms for binary two-stage robust optimization.
\newblock \emph{Computational Optimization and Applications}, 77\penalty0 (2):\penalty0 539--569, 2020.

\bibitem[Kasperski and Zieli{\'n}ski(2017)]{kasperski2017robust}
Adam Kasperski and Pawe{\l} Zieli{\'n}ski.
\newblock Robust recoverable and two-stage selection problems.
\newblock \emph{Discrete Applied Mathematics}, 233:\penalty0 52--64, 2017.

\bibitem[Lozano and Smith(2018)]{lozano2018binary}
Leonardo Lozano and J~Cole Smith.
\newblock A binary decision diagram based algorithm for solving a class of binary two-stage stochastic programs.
\newblock \emph{Mathematical Programming}, pages 1--24, 2018.

\bibitem[Lozano et~al.(2022)Lozano, Bergman, and Cire]{lozano2022constrained}
Leonardo Lozano, David Bergman, and Andre~A Cire.
\newblock Constrained shortest-path reformulations for discrete bilevel and robust optimization.
\newblock \emph{arXiv preprint arXiv:2206.12962}, 2022.

\bibitem[MacNeil and Bodur(2023)]{macneil2022leveraging}
Moira MacNeil and Merve Bodur.
\newblock Leveraging decision diagrams to solve two-stage stochastic programs with binary recourse and logical linking constraints.
\newblock \emph{European Journal of Operational Research}, 2023.

\bibitem[McElfresh et~al.(2019)McElfresh, Bidkhori, and Dickerson]{mcelfresh2019scalable}
Duncan~C McElfresh, Hoda Bidkhori, and John~P Dickerson.
\newblock Scalable robust kidney exchange.
\newblock In \emph{Proceedings of the AAAI Conference on Artificial Intelligence}, volume~33, pages 1077--1084, 2019.

\bibitem[Neumann(1928)]{neumann1928theorie}
J~v Neumann.
\newblock Zur theorie der gesellschaftsspiele.
\newblock \emph{Mathematische Annalen}, 100\penalty0 (1):\penalty0 295--320, 1928.

\bibitem[Postek and Hertog(2016)]{postek2016multistage}
Krzysztof Postek and Dick~den Hertog.
\newblock Multistage adjustable robust mixed-integer optimization via iterative splitting of the uncertainty set.
\newblock \emph{INFORMS Journal on Computing}, 28\penalty0 (3):\penalty0 553--574, 2016.

\bibitem[Rahmaniani et~al.(2017)Rahmaniani, Crainic, Gendreau, and Rei]{rahmaniani2017benders}
Ragheb Rahmaniani, Teodor~Gabriel Crainic, Michel Gendreau, and Walter Rei.
\newblock The {B}enders decomposition algorithm: A literature review.
\newblock \emph{European Journal of Operational Research}, 259\penalty0 (3):\penalty0 801--817, 2017.

\bibitem[Serra et~al.(2019)Serra, Raghunathan, Bergman, Hooker, and Kobori]{serra2019last}
Thiago Serra, Arvind~U Raghunathan, David Bergman, John Hooker, and Shingo Kobori.
\newblock Last-mile scheduling under uncertainty.
\newblock In \emph{International Conference on Integration of Constraint Programming, Artificial Intelligence, and Operations Research}, pages 519--528. Springer, 2019.

\bibitem[Subramanyam et~al.(2020)Subramanyam, Gounaris, and Wiesemann]{subramanyam2020k}
Anirudh Subramanyam, Chrysanthos~E Gounaris, and Wolfram Wiesemann.
\newblock K-adaptability in two-stage mixed-integer robust optimization.
\newblock \emph{Mathematical Programming Computation}, 12\penalty0 (2):\penalty0 193--224, 2020.

\bibitem[van Hoeve(2022)]{van2022graph}
Willem-Jan van Hoeve.
\newblock Graph coloring with decision diagrams.
\newblock \emph{Mathematical Programming}, 192\penalty0 (1):\penalty0 631--674, 2022.

\bibitem[Vayanos et~al.(2011)Vayanos, Kuhn, and Rustem]{vayanos2011decision}
Phebe Vayanos, Daniel Kuhn, and Ber{\c{c}} Rustem.
\newblock Decision rules for information discovery in multi-stage stochastic programming.
\newblock In \emph{2011 50th IEEE Conference on Decision and Control and European Control Conference}, pages 7368--7373. IEEE, 2011.

\bibitem[Yan and Kung(2018)]{yan2018robust}
Chiwei Yan and Jerry Kung.
\newblock Robust aircraft routing.
\newblock \emph{Transportation Science}, 52\penalty0 (1):\penalty0 118--133, 2018.

\bibitem[Yan{\i}ko{\u{g}}lu et~al.(2019)Yan{\i}ko{\u{g}}lu, Gorissen, and den Hertog]{yanikouglu2019survey}
{\.I}hsan Yan{\i}ko{\u{g}}lu, Bram~L Gorissen, and Dick den Hertog.
\newblock A survey of adjustable robust optimization.
\newblock \emph{European Journal of Operational Research}, 277\penalty0 (3):\penalty0 799--813, 2019.

\bibitem[Zeng and Zhao(2013)]{zeng2013solving}
Bo~Zeng and Long Zhao.
\newblock Solving two-stage robust optimization problems using a column-and-constraint generation method.
\newblock \emph{Operations Research Letters}, 41\penalty0 (5):\penalty0 457--461, 2013.

\end{thebibliography}

\newpage



\ECSwitch

\ECHead{Electronic Companion}

\section{Proofs}\label{appendix:proofs}

\paragraph{Proof of Lemma \ref{propn:conv_seperation}.}

We first show that $\conv(\mY \cap \mS(\bx)) \subseteq \conv(\mY) \cap \mS(\bx), \  \forall \bx \in \{0,1\}^m$. For any $\bhx \in \{0,1\}^m$, let $\bhy$ be an arbitrary point in $\conv(\mY \cap \mS(\bhx))$. By definition, $\bhy$ can be written as a convex combination of points $\by^1, \ldots, \by^R$ that are feasible in both $\mY$ and $\mS(\bhx)$. Since $\by^1, \ldots, \by^R \in \mY$, it must be true that $\bhy \in \conv(\mY)$. Similarly, since $\by^1, \ldots, \by^R \in \mS(\bhx)$, $\bhy \in \conv(\mS(\bhx)) = \mS(\bhx)$.

Next, we show that $\conv(\mY) \cap \mS(\bhx) \subseteq \conv(\mY \cap \mS(\bhx))$. We prove this by contradiction. Suppose there exists a vector $\bhy \in \conv(\mY) \cap \mS(\bhx)$ such that $\bhy \notin \conv(\mY \cap \mS(\bhx))$. Then, $\bhy$ must satisfy one of the two conditions:
\begin{itemize}
    \item The vector $\bhy$ is a binary vector, i.e., $\bhy \in \{0,1\}^n$. However, if $\bhy \in \{0,1\}^n$ and $\bhy \in \conv(\mY)$, then $\bhy \in \mY$. This implies that $\bhy \in \mY$ and $\bhy \in \mS(\bhx)$, which further implies that $\bhy \in \conv(\mY \cap \mS(\bhx))$. This is a contradiction. 
    \item The vector $\bhy$ must have at least one index $k$ where $\hat{y}_k \in (0,1)$. This implies that $\bhy$ must be a strict convex combination of a set of binary vectors $\by^1, \ldots, \by^r \in \mY$, that is, $\bhy = \sum_{r=1}^R \lambda_r \by^r$, $\sum_{r=1}^R \lambda_r = 1$, and $\lambda_r > 0 \ \forall r \in \{1, \ldots, R\}$. If all binary vectors $\by^1, \ldots, \by^r \in \mS(\bhx)$, then $\bhy \in \conv(\mS(\bhx)) = \mS(\bhx)$ which implies that $\bhy \in \conv(\mY \cap \mS(\bhx))$. This is a contradiction. Thus, it must be the case that there exists a binary vector $\by^\ell \in \{\by^1, \ldots, \by^r\}$ such that $\by^\ell \notin \mS(\bhx)$. Based on our definition of $\mS(\bhx)$, this implies that there is a violated constraint of the form $y^\ell_i \leq \hat{x}_j$, $y^\ell_{i} = \hat{x}_j$ or $y^\ell_i \geq \hat{x}_j$ for some $j \in \{ 1,\hdots,m\}$. Recall that $\by^\ell$ and $\hat{\bx}$ are binary vectors. Without loss of generality, suppose $y^\ell_i = 1$ which violates the constraint $y^\ell_i \leq \hat{x}_{j}$ when $\hat{x}_j = 0$. However, since $\bhy$ is a strict combination of binary vectors which includes $\by^\ell$, we have $\hat{y}_i > 0$, thus it must be true that $\hat{y}_i$ also violates this constraint. This implies that $\bhy \notin \mS(\bhx)$, which contradicts our initial assumption. This argument can be made for every value of $\hat{y}^\ell_i$ and $\hat{x}_j$ which violates one of the selectively adaptive constraints. This concludes our proof.  \Halmos
\end{itemize}

\medskip

\paragraph{Proof of Proposition \ref{thm1}.}
The main idea behind this proof comes from Proposition 1 of \cite{arslan2022decomposition}. We first note that Problem \eqref{model:1} is equivalent to
\begin{align*}  \underset{\bx \in \mX}{\text{min}} \ \underset{\bxi \in \Xi}{\text{max}} \  \underset{\by \in \conv(\mY \cap \mS(\bx))}{\text{min}} \quad \bc^\top \bx + \bxi^\top \by.
\end{align*}
Since both the maximization problem and the inner minimization problem are now over convex sets, we can apply the minimax theorem \citep{neumann1928theorie} to swap the order of these two operations and derive the equivalent reformulation of 
\begin{align*}
\underset{\bx \in \mX}{\text{min}} \ \underset{\by \in \conv(\mY \cap \mS(\bx))}{\text{min}} \ \underset{\bxi \in \Xi}{\text{max}} \quad \bc^\top \bx + \bxi^\top \by.
\end{align*}
Using Lemma \ref{propn:conv_seperation}, we have $\conv(\mY\cap \mS(\bx)) = \conv(\mY) \cap \mS(\bx)$, thus we can combine the two minimization operations into a single one that is solved with $\bx \in \mX, \, \by \in \conv(\mY)$, and $\by \in \mS(\bx)$. Finally, we can introduce the auxiliary variable $v$ to switch to the epigraph formulation of the maximization problem, which concludes the proof. \Halmos

\medskip

\paragraph{Proof of Proposition \ref{propn:constraint_decoupling}}

This proposition comes directly from the result that for any two sets $\mA, \mB \subseteq \mathbb{R}^n$,  it must be true that $\conv(\mA \cap \mB) \subseteq \conv(\mA) \cap \conv(\mB)$. As a quick proof, let $\bx$ be an arbitrary point in $\conv(\mA \cap \mB)$. By definition, $\bx$ can be written as a convex combination of feasible points in $\mA \cap \mB$. Since these feasible points are in both $\mA$ and $\mB$, they must also be in $\conv(\mA)$ and $\conv(\mB)$. Finally, since $\bx$ is a convex combination of these points, $\bx$ must also be in $\conv(\mA)$ and $\conv(\mB)$. This implies that $\conv(\mA \cap \mB) \subseteq \conv(\mA) \cap \conv(\mB)$. Furthermore, we also note that for any set $\mB'$ such that $\mB \subseteq \mB'$, it must be true that $\conv(\mB) \subseteq \conv(\mB')$, which implies that $\conv(\mA \cap \mB) \subseteq \conv(\mA) \cap \conv(\mB')$. 

By applying these two arguments, $\conv(\mY)$ must be a subset of the intersection of sets (i) $\rel(\mY)$, (ii) $\proj_\by(\NF(\mD^{\mJ_i})) \ \forall \mJ_i \in \{\mJ_1^1, \ldots, \mJ^k_1\}$, and (iii) $\proj_\by(\NF(\mD^{\mJ_i}_{\texttt{outer}})) \ \forall \mJ_i \in \{\mJ_1^2, \ldots, \mJ^2_q\}$. In other words, the intersection of these three sets is a valid outer approximation of $\conv(\mY)$. \Halmos

\newpage

\section{Details of the Capital Budgeting Problem}

The complete capital budgeting problem from \cite{arslan2022decomposition} is 
\begin{align*}
    \underset{(\bx, x_0) \in \mX}{\text{max}}  \underset{\xi \in \Xi}{\text{min}} \ \underset{(\by, y_0) \in \mY \cap \mS(\bx)}{\text{max}} \ \sum_{i \in \mN} \bigg( \sum_{j = 1}^M \frac{Q_{ij}\xi_j}{2} \bigg) \bar p_i ((1-f)x_i + fy_i) + \sum_{i \in \mN} \bar p_i((1-f)x_i + fy_i) - \gamma x_0 - \gamma \mu y_0
\end{align*}
where $\Xi := [-1,1]^M$ and where
\begin{align*}
\mY \cap \mS(\bx) = 
\left \{(\by, y_0, w_0) \in \{0,1\}^{N+2} \ \bigg | \ 
 \bc^\top \by \leq B + C_1w_0 + C_2y_0, \ w_0 = x_0, \ y_i \geq x_i \ \forall i \in \mN \right \}. 
\end{align*}

Note that $w_0$ is an auxiliary variable that is used to maintain the selective adaptability condition.

\subsection{Deriving the network flow formulation}

We show a step-by-step process to obtain our network flow formulation. First, recall that if we had a description for $\conv(\mY\cap\mS(\bx))$, then the problem can be rewritten as
\begin{align*}
    \underset{\substack{(\bx, x_0) \in \mX \\ (\by, y_0) \in \conv(\mY \cap \mS(\bx))}}{\text{max}} \ - \gamma x_0 - \gamma \mu y_0 + \sum_{i \in \mN} \bar p_i((1-f)x_i + fy_i) + \underset{\xi \in \Xi}{\text{min}} \  \sum_{i \in \mN} \bigg( \sum_{j = 1}^M \frac{Q_{ij}\xi_j}{2} \bigg) \bar p_i ((1-f)x_i + fy_i). 
\end{align*}

Given the uncertainty set $\Xi := [-1,1]^M$, the adversarial problem can be written as 
\begin{align*}
\underset{\bxi}{\text{min}} \quad &  \sum_{i \in \mN} \bigg( \sum_{j = 1}^M \frac{Q_{ij}\xi_j}{2} \bigg) \bar p_i ((1-f)x_i + fy_i) \\
\text{s.t.} \quad & \xi_i \geq -1,  \quad \quad \forall i \in \mM \\ 
& \xi_i \leq 1, \quad \quad \ \ \forall i \in \mM,
\end{align*}
where $\mM = \{1, \ldots, M\}$. The dual of this problem is 
\begin{align*}
\underset{\blambda^1, \blambda^2}{\text{max}} \quad &  -\sum_{i \in \mM} \lambda^1_m + \lambda^2_m \\
\text{s.t.} \quad & \lambda^1_m - \lambda^2_m \leq \sum_{i \in \mN} \bigg( \frac{Q_{i,m}}{2} \bigg) \bar p_i ((1-f)x_i + fy_i) \quad \forall m \in \mM \\
& \blambda^1, \blambda^2 \geq \bzero.
\end{align*}

Finally, combining this inner maximization problem with the outer maximization problem, we obtain our final network flow formulation of

\begin{align}
\begin{split}
\underset{(\bx, x_0), (\by,y_0), \bz, \blambda}{\text{max}} \quad & -\gamma x_0 + \sum_{i \in \mN} \bar p_i((1-f)x_i + fy_i) - (\sum_{m \in \mM} \lambda^1_m + \lambda^2_m) - \gamma\mu y_0\\
\text{s.t.} \quad & \lambda^1_m - \lambda^2_m = \sum_{i \in \mN} \bigg( \frac{Q_{i,m}}{2} \bigg) \bar p_i ((1-f)x_i + fy_i) \quad \forall m \in \mM,\\
& \bA \bz = \bb, \\ 
& y_i = \sum_{j \in l^+(i)} z_j, \quad \forall i \in \mN\\
& y_0 = \sum_{j \in l^+(N-1)} z_j, \\
& w_0 = \sum_{j \in l^+(N)} z_j,\\
& \by \geq \bx,\\
& w_0 = x_0,\\
& \bc^\top \bx \leq B + C_1x_0,\\
& (\bx, x_0) \in \{0,1\}^{N+1},\\
& \blambda^1, \blambda^2, \bz \geq \bzero.
\end{split}
\end{align}

\vspace{0.2cm}

\subsection{Deriving the $K$-adaptability formulation}\label{ECsub:KAdaptCapital}

We derive the $K$-adaptability formulation for the capital budgeting problem, which is based on Theorem 2 of \cite{hanasusanto2015k}. For a given solution $(\bx, \by^1, \ldots, \by^K)$, the adversarial problem can be written as
\begin{align*}
\underset{\bxi, \tau}{\text{min}} \quad &  \tau \\ 
\text{s.t.} \quad & \tau \geq - \gamma x_0 - \gamma \mu y^k_0 + \sum_{i \in \mN} \bigg( \sum_{j = 1}^M \frac{Q_{ij}\xi_j}{2} \bigg) \bar p_i ((1-f)x_i + fy^k_i) + \sum_{i \in \mN} \bar p_i((1-f)x_i + fy^k_i), \quad \forall k \in \mK \\
& \xi_i \geq -1,  \quad \quad \forall i \in \mM \\ 
& \xi_i \leq 1, \quad \quad \ \ \forall i \in \mM.
\end{align*}

The dual of this problem is 
\begin{align*}
\underset{\blambda^1, \blambda^2, \bpi}{\text{min}} \quad &  -  (\sum_{m \in \mM} \lambda^1_m + \lambda^2_m) + \sum_{k =1}^K \pi_k \bigg ( - \gamma x_0 - \gamma \mu y^k_0 + \sum_{i \in \mN} \bar p_i((1-f)x_i + fy^k_i) \bigg)\\ 
\text{s.t.} \quad & \lambda^1_m - \lambda^2_m = \sum_{k = 1}^K \sum_{i \in \mN} \bigg( \frac{Q_{i,m}}{2} \bigg) \bar p_i ((1-f)x_i \pi_k  + fy_i^k \pi_k ) \quad \forall m \in \mM \\
& \sum_{k = 1}^{K} \pi_k = 1\\ 
& \blambda^1, \blambda^2 \in \mathbb{R}^m, \bpi \in \mathbb{R}^k.
\end{align*}

Note that since $\sum_{k = 1}^K \pi_k = 1$, the model is equivalent to
\begin{align*}
\underset{\blambda^1, \blambda^2, \bpi}{\text{min}} \quad &  -  (\sum_{m \in \mM} \lambda^1_m + \lambda^2_m) - \gamma x_0 + (1-f)\sum_{i \in \mN} \bar p_ix_i + \sum_{k =1}^K \sum_{i \in \mN} f \bar p_i y^k_i \pi_k  - \gamma \mu y^k_0 \pi_k \\ 
\text{s.t.} \quad & \lambda^1_m - \lambda^2_m = \sum_{i \in \mN} \bigg(\frac{Q_{i,m}}{2}\bigg) \bar p_i (1-f) x_i + \sum_{k = 1}^K \sum_{i \in \mN} \bigg( \frac{Q_{i,m}}{2} \bigg) \bar p_i fy_i^k \pi_k  \quad \forall m \in \mM \\
& \sum_{k = 1}^{K} \pi_k = 1\\ 
& \blambda^1, \blambda^2 \in \mathbb{R}^m_+, \bpi \in \mathbb{R}^k_+.
\end{align*}

Finally, because of the bilinearity, we can replace every instance of $\pi_ky^k_i$ with $q^k_i$. We thus arrive at the final $K$-adaptability formulation of 
\begin{alignat*}{2}
\text{min} \quad &  -  (\sum_{m \in \mM} \lambda^1_m + \lambda^2_m) - \gamma x_0 + (1-f)\sum_{i \in \mN} \bar p_ix_i + \sum_{k =1}^K \sum_{i \in \mN} f \bar p_i q^k_i  - \gamma \mu q^k_0\\ 
\text{s.t.} \quad & \lambda^1_m - \lambda^2_m = \sum_{i \in \mN} \bigg(\frac{Q_{i,m}}{2}\bigg) \bar p_i (1-f) x_i + \sum_{k = 1}^K \sum_{i \in \mN} \bigg( \frac{Q_{i,m}}{2} \bigg) \bar p_i fq_i^k, \quad && \forall m \in \mM \\
& q^k_i \leq y^k_i, q^k_i \leq \pi_k, q^k_i \geq \pi_k + y^k_i - 1, \quad && \forall i \in \{1, \ldots, |\mS|\}, k \in \{1, \ldots, K\} \\
&  \bc^\top \by^k \leq B + C_1x_0 + C_2y^k_0, \quad && \forall k \in \{1, \ldots, K\} \\
& \by^k \geq \bx, \quad && \forall k \in \{1, \ldots, K\}\\
& \bc^\top \bx \leq B + C_1x_0 \\
& \sum_{k = 1}^{K} \pi_k = 1\\ 
& \blambda^1, \blambda^2 \in \mathbb{R}^m_+, \ \bpi \in \mathbb{R}^k_+, \ \bq^1, \ldots, \bq^k \in \mathbb{R}^n_+ \\ 
& \bx, \by^1, \ldots, \by^k \in \{0,1\}^n.
\end{alignat*}

\vspace{0.1cm}

\subsection{A brief comparison with branch-and-price results}\label{ECsub:BPresults}

We briefly compare the results of our network flow models against the solution times presented in \cite{arslan2022decomposition} for the capital budgeting instances. As discussed in the main body of the paper, branch-and-price algorithms represent a family of methods that can be effective if implemented properly, but may require significant implementation effort (e.g., fine-tuning of many intermediate steps, see \cite{cire2019network, de2022arc}). The main purpose of our comparison is simply to show that our approach, which is model-based and can be solved directly using any standard solver, is a viable and promising alternative for solving these problems.

\begin{table}[h]
\centering
\setlength{\tabcolsep}{12pt}
{\renewcommand{\arraystretch}{1.2}%
\begin{tabular}{c|llll|r}
 \toprule
Instances & $Q = 0$ & $Q = 1$ & $Q = 3$ & $Q = 5$ & \ \ \  B\&P  \\
\midrule 
 n = 10  & $0.1$s & - & - & - & $0.1$s \\ 
 n = 20  & $4$s & - & - & - & $1.8$s \\ 
 n = 30  & $104$s & $16$s $(0.3\%)$ & - & - & $>268$s \\ 
 n = 40  & $589$s & $78$s $(0.1\%)$ & $12$s $(0.4\%)$ & - & $>451$s  \\ 
 n = 50  & - & $150$s $(0.08\%^*)$ & $30$s $(0.2\%^*)$ & $12$s  $(0.3\%^*)$ & $116$s \\
\bottomrule
\end{tabular}}
\vspace{0.1cm}
\caption{Solution time of the network flow models compared to the branch-and-price algorithm. Notes: Solution times of B\&P algorithm are taken from \cite{arslan2022decomposition}. Two of the $n = 30$ instances and five of the $n = 40$ instances could not be solved within the one-hour time limit using the B\&P algorithm. For the $n = 30$ and $n = 40$ instances, we give the true optimality gap of the solutions in the brackets. For the $n = 50$ instances, we present the model-based optimality gap (denoted using $*$), which is an upper bound on the true optimality gap.}\label{table:vs_BP}
\end{table}

Table \ref{table:vs_BP} highlights the total time to compile and solve the network flow models, as well as the solution time of the branch-and-price (B\&P) algorithm reported in \cite{arslan2022decomposition}. For the network flow models, the times are taken directly from Table \ref{table:RDD_times} in the main body of the paper, and we censor some values here for ease of discussion. The values in the brackets next to the solution times denote the average optimality gap of the solutions, which are taken directly from Table \ref{table:RDD_optgap}. We remind the reader that there are 60 instances for each value of $n = 10, 20, 30, 40, 50$.

First, we observe that for the $n = 10$ and $n = 20$ instances, both methods can solve the instances exactly within a few seconds. For the $n = 30$ instances, our exact network flow model takes an average of 104 seconds to solve while the B\&P algorithm takes an average of 268 seconds with 2 instances that cannot be solved within the one-hour time limit imposed in the paper. For the instances where $n = 40$ and $n = 50$, our exact formulations become more difficult to compile and solve, and their solution times often exceed those of the B\&P algorithm. However, the approximation techniques in our paper provide a method for drastically reducing solution times while sacrificing little in solution quality. For example, for the $n = 40$ instances, we can generate significantly smaller network flow models with $Q = 1$ ($Q = 3$) that solves in an average of 78 (12) seconds while generating solutions that have an average optimality gap of 0.1\% (0.4\%). In comparison, the B\&P algorithm requires an average of 451 seconds with 5 instances that cannot be solved within the one-hour time limit. Similar observations can be made for the $n = 50$ instances. 

In summary, we remark that both exact approaches will become more difficult as the size of instances increase (this will naturally be the case for any exact approach for solving ARBO problems). We observe that the solution times of our approach are comparable to those of branch-and-price, but more importantly, we have proposed a rigorous framework for generating approximate formulations that are significantly faster to solve while sacrificing little in terms of solution quality.

\newpage

\section{Details of the Robust Assignment Problem}\label{app:RAP}

\subsection{Recursive formulation}\label{appsub:RAP_recursion}

In order to generate the exact decision diagram for the set of constraints
\begin{alignat*}{2}
& \sum_{\ell \in \mL(m)} a_\ell y_{\ell,m} \leq b_m, \quad && \forall m \in \{1, \ldots, M\} \\[0.1cm]
& \sum_{m \in \mM(\ell)} y_{\ell,m} \leq 1, \quad && \forall \ell \in \{1, \ldots, L\} \\[0.1cm]
& \ \by \in \{0,1\}^{|\mS|},
\end{alignat*}
we can extend Example \ref{ex:BDD_construction} to consider state vectors. For simplicity, we abuse notation slightly and consider decision vector $\by$ in the form of $\by = [y_1, \ldots, y_{|\mS|}]$, where for some element $y_{k=(\ell,m)}$, $agent(y_{k})$ is a zero vector with a single entry $a_\ell y_{k}$ at the $m$-th element, where $task(y_k)$ is a zero vector with a single entry $y_{k}$ at the $\ell$-th element. The recursive formulation can then be given with initial state $\bS_1 = \bzero$, state-transition function $T_i(\bS, y_k) = \bS + [agent(y_k); task(y_k)]$ and feasible action space $\mQ_i(\bS) = \{y_k \in \{0,1\} | T_i(\bS, y_k) \leq [\bb ; \bone] \}$.

\subsection{Complete network flow formulation}

Recall that the robust assignment problem considered is a max-min-max problem. Through a simple transformation, we can consider an equivalent min-max-min formulation (i.e., by applying a negative in the objective function and taking the negative of the optimal objective value). We use this min-max-min model to formulate our network flow problem. We show the step-by-step process below.

Recall that the uncertainty set to the robust assignment problem is defined as
\begin{align*}
\Xi = \left\{\bxi \; \bigg | \; \sum_{i = 1}^{|\mS|} \big |\xi_i/\xi^0_i - 1 \big| \leq 0.1|\mS|, \ \big |\xi_i/\xi^0_i - 1 \big | \leq 0.5 \ \ \forall i \in \{1, \ldots, |\mS|\} \right\}.  
\end{align*}

Given this uncertainty set, we can then write the adversarial problem as 
\begin{alignat*}{2}
\underset{\bxi, \bw}{\text{max}} \quad & -\bxi^\top \by\\
\text{s.t.} \quad & w_i \geq \xi_i/\xi^0_i - 1, \quad \ && \forall i \in \{1, \ldots, |\mS|\}\\
& w_i \geq 1 - \xi_i/\xi^0_i, \quad && \forall i \in \{1, \ldots, |\mS|\}\\
& w_i \leq 0.5, \quad && \forall i \in \{1, \ldots, |\mS|\}\\
& \sum_{i=1}^{|\mS|} w_i \leq 0.1|\mS|.
\end{alignat*}
The dual of this problem is
\begin{alignat*}{2}
\underset{}{\text{min}} \quad & 0.1|\mS|\alpha + \sum_{i=1}^{|\mS|} (\lambda_i^1 - \lambda_i^2 + 0.5\pi_i)\\
\text{s.t.} \quad & \lambda^1_i - \lambda^2_i = -\xi^0_iy_i, \quad && \hspace*{-1cm} \forall i \in \{1, \ldots, |\mS|\}\\
& \pi_i - \lambda_i^1 - \lambda^2_i + \alpha = 0, \quad && \hspace*{-1cm} \forall i \in \{1, \ldots, |\mS|\}\\
& \blambda^1, \blambda^2, \bpi \in \mathbb{R}^{|\mS|}_+ , \  \alpha \geq 0.
\end{alignat*}
Finally, the complete network flow formulation of the adaptive robust optimization problem is
\begin{alignat*}{2}
\underset{}{\text{min}} \quad & 0.1|\mS|\alpha + \sum_{i=1}^{|\mS|} (\lambda_i^1 - \lambda_i^2 + 0.5\pi_i)\\
\text{s.t.} \quad &  \lambda^1_i - \lambda^2_i = -\xi^0_iy_i, && \hspace*{-2cm} \forall i \in \{1, \ldots, |\mS|\}\\
& \pi_i - \lambda_i^1 - \lambda^2_i + \alpha = 0,  && \hspace*{-2cm} \forall i \in \{1, \ldots, |\mS|\}\\
& \by \in \proj(\text{NF}(\mD)) \\ 
& \by \leq \bx \\
&  ||\bx||_1 \leq \beta \\
& \bx \in \{0,1\}^{|\mS|}, \  \by, \blambda^1, \blambda^2, \bpi \in \mathbb{R}^{|\mS|}_+ , \  \alpha \geq 0.
\end{alignat*}

\subsection{K-adaptability formulation}\label{appsub:RAP_Kadapt}

We similarly outline the step-by-step process for formulating the $K$-adaptability model, which is based on \cite{hanasusanto2015k}. First, for a given solution $(\bx, \by^1, \ldots, \by^K)$, the epigraph formulation of the adversarial maximization problem is
\begin{alignat*}{2}
\underset{\bxi, \bw, \tau}{\text{max}} \quad & \tau\\
\text{s.t.} \quad & \tau \leq -\bxi^\top \by^k \quad && \forall k \in \{1,\ldots, K\}\\
& w_i \geq \xi_i/\xi^0_i - 1, \quad && \forall i \in \{1, \ldots, |\mS|\}\\
& w_i \geq 1 - \xi_i/\xi^0_i, \quad && \forall i \in \{1, \ldots, |\mS|\}\\
& w_i \leq 0.5, \quad && \forall i \in \{1, \ldots, |\mS|\}\\
& \sum_{i=1}^{|\mS|} w_i \leq 0.1|\mS|.
\end{alignat*}

This problem has a dual of
\begin{alignat*}{2}
\underset{}{\text{min}} \quad & 0.1|\mS|\alpha + \sum_{i=1}^{|\mS|} (\lambda_i^1 - \lambda_i^2 + 0.5\pi_i)\\
\text{s.t.} \quad & \lambda^1_i - \lambda^2_i = - \sum_{k = 1}^K \xi^0_i\mu_ky^k_i, \quad && \forall i \in \{1, \ldots, |\mS|\}\\
& \pi_i - \lambda_i^1 - \lambda^2_i + \alpha = 0, \quad && \forall i \in \{1, \ldots, |\mS|\}\\
& \sum_{k = 1}^K \mu_k = 1 \\
& \blambda^1, \blambda^2, \bpi \in \mathbb{R}^{|\mS|}_+ , \  \bmu \in \mathbb{R}^{K}, \  \alpha \geq 0.
\end{alignat*}
Note that there are bilinear terms of the form $\mu_ky^k_i$. However, since $\by^k \in \{0,1\}^n$, we can replace each bilinear term with auxiliary variable $q^k_i$ and add constraints
\begin{align*}
q^k_i \leq y^k_i, q^k_i \leq \mu_k, q^k_i \geq \mu^k + y^k_i - 1, \quad \forall i \in \{1, \ldots, |\mS|\}, k \in \{1, \ldots, K\}.
\end{align*}
Since the outer problem of choosing $\bx, \by^1, \ldots, \by^K$ is also a minimization problem, we obtain the following $K$-adaptability formulation
\begin{alignat*}{2}
\underset{}{\text{min}} \quad & 0.1|\mS|\alpha + \sum_{i=1}^{|\mS|} (\lambda_i^1 - \lambda_i^2 + 0.5\pi_i)\\
\text{s.t.} \quad & \sum_{k = 1}^K \mu_k = 1\\
&  \lambda^1_i - \lambda^2_i = - \sum_{k = 1}^K \xi^0_iq^k_i, \quad && \forall i \in \{1, \ldots, |\mS|\}\\
& \pi_i - \lambda_i^1 - \lambda^2_i + \alpha = 0, \quad && \forall i \in \{1, \ldots, |\mS|\}\\
& q^k_i \leq y^k_i, q^k_i \leq \mu_k, q^k_i \geq \mu_k + y^k_i - 1, \quad && \forall i \in \{1, \ldots, |\mS|\}, \ k \in \{1, \ldots, K\}\\
& \by^k \leq \bx, \quad && \forall k \in \{1, \ldots, K\}\\
& ||\bx||_1 \leq \beta\\
& \by^1, \ldots, \by^K \text{ are feasible assignments}\\
& \bx, \by^1, \ldots, \by^K, \in \{0,1\}^{|\mS|} \\   & \bq^1, \ldots, \bq^K, \blambda^1, \blambda^2, \bpi \in \mathbb{R}^{|\mS|}_+ , \ \bmu \in \mathbb{R}^K, \  \alpha \geq 0.
\end{alignat*}

\vspace{0.2cm}

\end{document}